\RequirePackage{fix-cm}
\documentclass[11pt]{article}

\usepackage[a4paper,margin=1in]{geometry}
\usepackage{fontspec}
\usepackage{polyglossia}
\setmainlanguage{english}
\setotherlanguage{russian}
\defaultfontfeatures{Ligatures=TeX}
\newfontfamily\russianfont{Tempora-Regular.otf}[
  BoldFont=Tempora-Bold.otf,
  ItalicFont=Tempora-Italic.otf,
  BoldItalicFont=Tempora-BoldItalic.otf
]
\newfontfamily\cyrillicfont{Tempora-Regular.otf}[
  BoldFont=Tempora-Bold.otf,
  ItalicFont=Tempora-Italic.otf,
  BoldItalicFont=Tempora-BoldItalic.otf
]
\newfontfamily\cyrillicfonttt{FreeMono.otf}[
  BoldFont=FreeMonoBold.otf,
  ItalicFont=FreeMonoOblique.otf,
  BoldItalicFont=FreeMonoBoldOblique.otf
]
\usepackage{amsmath,amssymb,amsthm,mathtools,bm,mathrsfs}
\input{ursfs.fd}
\DeclareFontShape{U}{rsfs}{m}{n}{<-> rsfs10}{}
\DeclareMathSizes{5}{3.5}{2.8}{2.5}
\DeclareMathSizes{6}{4.5}{3.5}{3}
\DeclareMathSizes{7}{5.5}{4}{3.5}
\DeclareMathSizes{8}{6.5}{4.875}{4.06}
\DeclareMathSizes{9}{7.5}{5.83}{4.17}
\DeclareMathSizes{10}{8.5}{5.95}{4.25}
\DeclareMathSizes{10.95}{9.45}{6.9}{5.18}
\DeclareMathSizes{12}{10.5}{7}{5.25}
\DeclareMathSizes{12.5}{11}{7.7}{5.5}
\DeclareMathSizes{14.4}{12.9}{8.96}{6.27}
\DeclareMathSizes{16}{14.5}{10.15}{7.25}
\DeclareMathSizes{17.28}{15.78}{10.96}{8.22}
\DeclareMathSizes{20.74}{19.24}{13.36}{9.28}
\DeclareMathSizes{24.88}{23.38}{16.24}{11.28}
\usepackage{etoolbox}
\makeatletter
\patchcmd{\text@}{\f@size}{\tf@size}{}{\PackageError{mathsize}{Cannot resize display math text}{}}
\patchcmd{\text@}{\f@size}{\tf@size}{}{\PackageError{mathsize}{Cannot resize inline math text}{}}
\makeatother
\usepackage{booktabs,array,enumitem,microtype,needspace,graphicx,float}
\usepackage{xcolor}
\usepackage{titlesec}
\usepackage{tikz}
\usetikzlibrary{arrows.meta,decorations.markings}
\usepackage[numbers,sort&compress]{natbib}
\usepackage{xurl}
\usepackage[hypertexnames=false]{hyperref}
\definecolor{ink}{HTML}{000000}
\definecolor{softgray}{HTML}{F3F4F6}
\definecolor{rulegray}{HTML}{9CA3AF}
\definecolor{geoPtwo}{HTML}{17365D}
\definecolor{geoPthreehalf}{HTML}{2F5F8F}
\definecolor{geoPsixfifths}{HTML}{6F8FB5}
\definecolor{geoPone}{HTML}{9CBAD0}

\tikzset{
  geo axis/.style={->,draw=black!70,line width=.45pt},
  geo attainable/.style={solid,line width=.9pt,line cap=round,line join=round},
  geo unattainable/.style={draw=black!75,densely dashed,line width=.8pt,
    line cap=round,line join=round},
  geo support/.style={draw=black!48,line width=.55pt},
  geo excluded/.style={circle,draw=black,fill=white,line width=.65pt,
    inner sep=1.25pt},
  geo included/.style={circle,draw=black,fill=black,inner sep=1.05pt}
}

\hypersetup{
  colorlinks=true,
  allcolors=black,
  pdftitle={Exponential Convergence of the Lifted Gradient Flow of the Bass Functional},
  pdfauthor={Danila Litvinov},
  pdfsubject={Complete English version followed by the Russian version},
  pdfkeywords={Bass functional, gradient flow, martingale optimal transport, entropy, strong convexity}
}

\titleformat{\section}
  {\normalfont\bfseries\color{ink}\hyphenpenalty=10000\exhyphenpenalty=10000\fontsize{16pt}{19pt}\selectfont}
  {\thesection}{0.75em}{}
\titleformat{\subsection}
  {\normalfont\bfseries\color{ink}\hyphenpenalty=10000\exhyphenpenalty=10000\fontsize{12.5pt}{15pt}\selectfont}
  {\thesubsection}{0.65em}{}
\titlespacing*{\section}{0pt}{3.5ex plus 1ex minus .2ex}{1.4ex}
\titlespacing*{\subsection}{0pt}{2.4ex plus .7ex minus .2ex}{.8ex}

\newtheoremstyle{basslemma}{6pt}{6pt}{\normalfont}{}{\itshape}{.}{0.5em}%
{\thmname{#1}\thmnumber{ #2}\thmnote{ (#3)}}
\newtheoremstyle{bassremark}{6pt}{6pt}{\normalfont}{}{\bfseries}{.}{0.5em}{}
\theoremstyle{plain}
\newtheorem{theoremen}{Theorem}[section]
\newtheorem{propositionen}[theoremen]{Proposition}

\newtheorem{theoremru}{Теорема}[section]
\newtheorem{propositionru}[theoremru]{Утверждение}

\theoremstyle{basslemma}
\newtheorem{lemmaen}[theoremen]{Lemma}
\newtheorem{lemmaru}[theoremru]{Лемма}
\theoremstyle{bassremark}
\newtheorem{remarken}[theoremen]{Remark}
\newtheorem{remarkru}[theoremru]{Замечание}

\newcommand{\E}{\mathbb{E}}
\newcommand{\G}{\mathcal{G}}
\newcommand{\norm}[1]{\left\lVert #1\right\rVert}
\newcommand{\ip}[2]{\left\langle #1,#2\right\rangle}
\newcommand{\R}{\mathbb{R}}

\DeclareMathOperator*{\supp}{supp}

\makeatletter
\renewenvironment{proof}[1][\(\blacktriangleright\)]{%
  \par
  \pushQED{\qed}%
  \normalfont \topsep6\p@\@plus6\p@\relax
  \trivlist
  \item[\hskip\labelsep #1]\ignorespaces
}{%
  \popQED\endtrivlist\@endpefalse
}
\makeatother

\setlist[itemize]{leftmargin=1.5em,itemsep=.35em,topsep=.45em}
\allowdisplaybreaks[2]

\newcommand{\makearticletitle}[2]{%
  \newpage\null\vskip 2em%
  \begin{center}%
    {\LARGE #1\par}%
    \vskip 1.5em%
    {\large #2\par}%
  \end{center}%
  \par\vskip 1.5em%
}
\newcommand{\makerussiantitle}{%
  \makearticletitle{\vspace{-1.2cm}\textbf{Экспоненциальная сходимость\\[2mm]
  поднятого градиентного потока\\[2mm] функционала Басса}}{Данила Литвинов}%
}
\newcommand{\Law}{\mathscr L}
\newcommand{\MCov}{\operatorname{MCov}}
\newcommand{\cH}{\mathscr H}

\usepackage{longtable}
\usepackage{tabularx}
\usetikzlibrary{positioning,calc,fit,backgrounds}

\newcommand{\makeenglishtitle}{%
  \makearticletitle{\vspace{-1.2cm}\textbf{Exponential Convergence\\[2mm]
  of the Lifted Gradient Flow\\[2mm] of the Bass Functional}}{Danila Litvinov}%
  \pdfbookmark[0]{English version}{language-en}%
}
\let\originalmakerussiantitle\makerussiantitle
\renewcommand{\makerussiantitle}{%
  \originalmakerussiantitle
  \pdfbookmark[0]{Русская версия}{language-ru}%
}

\begin{document}
\selectlanguage{english}
\color{ink}
\makeenglishtitle
\vspace{-1.2em}

\begin{abstract}
A gradient flow may converge to a minimizer without a known rate of convergence.
We study this question for the lifted Bass functional in martingale optimal
transport. Suppose that the source and target measures are in convex order
and form an irreducible pair. Assume that the target support is compact and
convex, and that its density is bounded above and bounded away from zero
almost everywhere on the interior of the support. Suppose also that the
source support lies strictly inside the target support. Under these
assumptions, we prove multidimensional exponential estimates for bounded
initial states. A single positive constant controls, at every nonnegative
time, both the decay of the $L^2$ distance to a global minimizer and the decay
of the functional gap. In both estimates, the prefactor is exactly the
corresponding initial error.

The main step is strong convexity on sets with a common $L^\infty$ bound
and a fixed mean. We
obtain it by introducing an entropic problem and interpolating the optimal
potentials at two endpoints. The second derivative of the comparison
function is expressed as a variance. A two-copy argument eliminates the
unknown change in the potential, while the interior geometry of the target
and Gaussian averaging yield a positive bound uniform in the entropic
parameter. This bound passes to the original functional as the entropic
regularization is removed. A known qualitative theorem places the entire trajectory in a
single bounded set, allowing us to obtain both rates.

The result extends the exponential-convergence theorem of Backhoff-Veraguas,
Pammer, and Schachermayer to dimensions $d\ge2$. No smoothness of the density
or the target boundary is required. The final result is accompanied by a
formal verification in Lean~4; the use of AI is described in
Section~\ref{en:sec:ai}.
\end{abstract}



\section{Introduction}\label{en:sec:intro}

A gradient flow decreases the value of a functional. This decrease alone,
however, does not say how quickly the trajectory approaches a minimizer.
Even after convergence has been established, a separate question remains:
can the error be bounded in terms of time and the initial state?

We consider this question for the lifted Bass functional arising in
martingale optimal transport. Backhoff-Veraguas, Pammer, and
Schachermayer~\citep{en:bps2025} proved strong $L^2$ convergence of this flow
in several dimensions under suitable assumptions on the underlying
measures. They obtained an exponential rate in one dimension. We prove
multidimensional exponential estimates under the original geometric and
probabilistic assumptions.

Why does this problem arise in martingale transport? The initial and final
distributions of a stochastic process are prescribed. The martingale
condition requires the conditional mean of a future state to equal the
current state. Even when the distributions are fixed, many admissible
couplings may exist, so a principle is needed to select a particular
construction. For a Bass martingale, finding a suitable source law becomes
a variational problem, and the lifted functional allows it to be studied
in a Hilbert space~\citep{en:bps2025}. This is why the quantitative behavior
of its gradient descent is of interest in its own right.

\subsection{The flow under consideration}

Let $X$ have the prescribed source law $\mu$, and let $\Gamma$ be an
independent standard Gaussian vector. The variable in the problem is a
random vector $Z$ measurable with respect to $\G=\sigma(X)$, viewed as an
element of $\cH=L^2(\G;\R^d)$. Gaussian noise is added to $Z$, after which
a random vector $Y$ with the prescribed target law $\nu$ is chosen.
The lifted functional maximizes the expected inner product of $Z+\Gamma$
and $Y$, while retaining the dependence between $Z$ and the original $X$
in a separate term:
\[
V(Z)=
\sup_{\substack{Y\in L^2(\sigma(X,\Gamma);\R^d)\\\Law(Y)=\nu}}
\left\{\E\ip{Z+\Gamma}{Y}-\E\ip{Z}{X}\right\}.
\]
Thus, both the law of $Z$ and its dependence on $X$ matter. Replacing the
second term by the maximal covariance between their laws would, in
general, change the problem.

We study the continuous-time equation
\[
\dot Z_t=-\nabla V(Z_t),\qquad Z_{t=0}=Z_0.
\]
Here $t$ is the time parameter of the motion toward a minimizer of the
functional. It differs from the time parameter of the martingale between
the initial and final distributions. The full setting, the gradient
formula, and all assumptions on the measures are given in the next section.

\subsection{From convergence to a rate}

Let $Z^\star$ be the limit of the trajectory. We want to control two errors:
the distance to this limit and the remaining functional gap. The aim is to
find $\kappa>0$ such that, simultaneously for every $t\ge0$,
\[
\begin{aligned}
\norm{Z_t-Z^\star}_2
&\le\norm{Z_0-Z^\star}_2e^{-\kappa t},\\
V(Z_t)-\inf_{\cH}V
&\le\bigl(V(Z_0)-\inf_{\cH}V\bigr)e^{-\kappa t}.
\end{aligned}
\]
In both estimates, the prefactor is exactly the corresponding initial
error. An estimate with an arbitrary additional prefactor, or one valid
only for sufficiently large times, would answer a weaker question.

We obtain these estimates for an irreducible pair of measures when the
target support is compact and convex, its density is bounded above and
bounded away from zero almost everywhere on the interior of that support,
and the source support lies strictly inside the target support. The
initial state is assumed to be bounded. The constant $\kappa$ may depend
on the underlying data and the initial state. We do not claim that it is
optimal or independent of all the data.

The mean constraint also matters. Convex order implies that $\mu$ and
$\nu$ have the same barycenter. Adding a constant vector to $Z$ therefore
leaves $V$ unchanged. Uniqueness of the minimizer must be understood within
the class of states with a prescribed mean. The minimum value itself,
however, is taken over the whole space $\cH$, including unbounded random
vectors.

\subsection{The main idea of the proof}

To obtain an exponential rate, it is natural to seek strong convexity.
In our problem, the required estimate holds on sets of the form
\[
\{Z\in\cH:\ \norm{Z}_\infty\le R,\ \E Z=\bar z\}.
\]
For each fixed radius $R$, we will obtain a positive curvature constant.
The known qualitative theorem provides a single radius that bounds the
entire trajectory. The same curvature constant can then be used at every
time.

The main difficulty is to establish this quantitative curvature. We
introduce entropy temporarily and take the optimal entropic potentials
at two source states. By linearly interpolating these endpoint potentials,
we construct an upper comparison for the value of the functional. Their
values agree at the endpoints; the intermediate potential need not be
optimal.

The second derivative of the comparison is expressed as a variance.
It contains an unknown difference of target potentials. Comparing two
independent copies eliminates this term. The remaining quantity is
related to the Hessian of a logarithmic integral. The geometry of the
target support and the lower bound on its density yield a positive
integrated bound for this Hessian. Gaussian averaging then transfers the
bound to the original space of random vectors.

\begin{figure}[H]
\centering
\begin{tikzpicture}[
  font=\small,
  bass box/.style={draw=black!45,rounded corners=2pt,line width=.5pt,
    align=center,inner sep=7pt,minimum height=1.02cm},
  bass step/.style={bass box,text width=4.7cm,fill=geoPtwo!6},
  bass known/.style={bass box,text width=4.7cm,fill=black!3},
  bass arrow/.style={-{Stealth[length=2mm]},draw=black!65,line width=.65pt}
]
\node[bass box,text width=10.9cm,fill=black!2] (data) at (0,2.5)
  {The given pair of measures and the lifted functional \(V\)};
\node[bass step] (entropy) at (-3.05,.9)
  {Entropic comparison\\and two independent copies};
\node[bass known] (qual) at (3.05,.9)
  {Known result:\\flow, boundedness,\\strong convergence};
\node[bass step] (curvature) at (-3.05,-.85)
  {Positive curvature \(m_R\)\\on \(\norm{Z}_\infty\le R\), \(\E Z=\bar z\)};
\node[bass known] (radius) at (3.05,-.85)
  {One radius for the entire orbit:\\\(\norm{Z_t}_\infty\le R\), \(t\ge0\)};
\node[bass box,text width=10.9cm,fill=geoPtwo!9] (rate) at (0,-2.8)
  {Strong monotonicity along the orbit\\
   The \(L^2\) error and the functional gap decay exponentially};
\draw[bass arrow] (data.south) -- ++(0,-.24) -| (entropy.north);
\draw[bass arrow] (data.south) -- ++(0,-.24) -| (qual.north);
\draw[bass arrow] (entropy) -- (curvature);
\draw[bass arrow] (qual) -- (radius);
\draw[bass arrow] (curvature.south) -- ++(0,-.32) -| ([xshift=-2.4cm]rate.north);
\draw[bass arrow] (radius.south) -- ++(0,-.32) -| ([xshift=2.4cm]rate.north);
\end{tikzpicture}
\caption{Structure of the rate proof. The quantitative argument of this paper is
on the left; the known qualitative result used in the proof is on the right.
The arrows indicate dependencies: the curvature estimate is applied with
a radius that bounds the entire orbit.}
\label{en:fig:bass-route}
\end{figure}
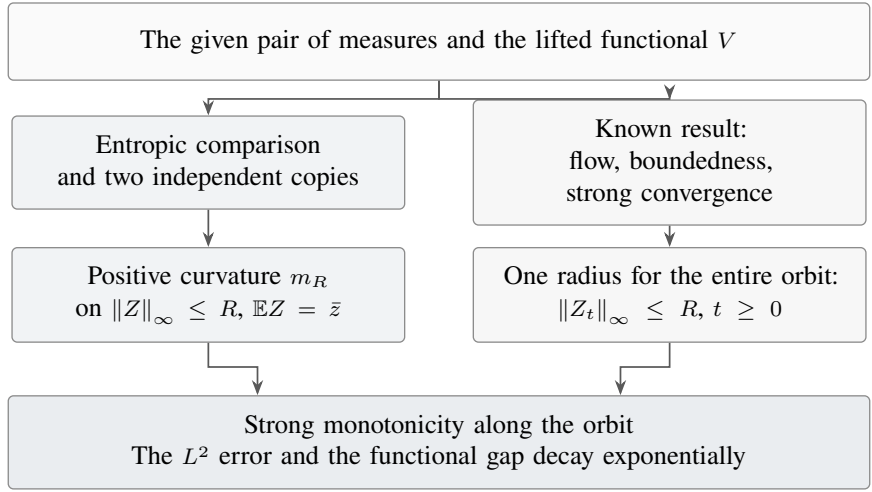

The resulting constant is independent of the entropic parameter. It
therefore survives as the entropy vanishes. Strong monotonicity of the
gradient along the bounded trajectory then gives both required estimates.
The diagram in Figure~\ref{en:fig:bass-route} shows where the known
qualitative convergence is used and where the new quantitative step enters.

\subsection{Organization of the paper}

Section~2 states the problem and the main result precisely and identifies
the properties of the known flow that we use. Section~3 introduces the
entropic comparison and justifies the passage to the original functional.
Section~4 establishes quantitative curvature, including the two-copy
estimate. Section~5 derives exponential convergence. Section~6 discusses
the significance of the result and its relation to other work. Section~7
describes the formal verification and how to reproduce the source checks;
Section~8 explains the use of AI. Section~9 concludes the paper.

\section{The setting and the main result}\label{en:sec:setting}

\subsection{The fixed data}

Let $d\ge2$, and let $\mu,\nu\in\mathcal P_2(\R^d)$ be probability measures
with finite second moments. We impose the following assumptions.
\begin{enumerate}[label=(A\arabic*),leftmargin=2.7em]
\item $\mu\le_{cx}\nu$: for every convex function for which the corresponding
integrals are defined, the integral against $\mu$ does not exceed the integral
against $\nu$.
\item The pair $(\mu,\nu)$ is irreducible: for any measurable sets
$A,B\subset\R^d$ with $\mu(A)>0$ and $\nu(B)>0$, there is a martingale coupling
$\pi$ of these measures such that $\pi(A\times B)>0$.
The martingale condition means that $\int y\,\pi_x(dy)=x$ for $\mu$-almost every
$x$. The coupling in this condition may depend on both sets.
\item The set $C=\supp\nu$ is compact and convex, and
$\supp\mu\subset I:=\operatorname{int}C$.
\item The measure $\nu$ is absolutely continuous with respect to
$d$-dimensional Lebesgue measure. Its density on $I$ satisfies
$\underline f\le f_\nu\le\overline f$ almost everywhere for some constants
$0<\underline f\le\overline f<\infty$.
\end{enumerate}

The last assumption in particular forces $C$ to be full-dimensional.
The boundary of a convex body has zero Lebesgue measure, so $\nu(I)=1$.
Since $\supp\mu$ is compact and lies in the interior of $C$, it is separated
from the boundary by a positive distance. We assume no smoothness of the
boundary or the density. These are the assumptions for the multidimensional
extension of Theorem~1.9 in~\citep{en:bps2025}. Here irreducibility is stated in
terms of couplings. To pass from a formulation in terms of a martingale,
one takes its joint law at times $0$ and $1$. Conversely, a coupling can be
realized as a two-step martingale with the same endpoints.

On an arbitrary probability space, fix independent random vectors
$X\sim\mu$ and $\Gamma\sim\gamma=N(0,I_d)$. Set
\[
 \G=\sigma(X),\qquad \mathcal F=\sigma(X,\Gamma),\qquad
 \cH=L^2(\G;\R^d).
\]
Equalities between random variables are understood up to almost sure
equality. The inner product on $\cH$ is
$\ip{Z}{W}_{\cH}=\E\ip{Z}{W}$; we write $\norm{Z}_2$ for the associated norm.

\subsection{Why the lifted functional is needed}

For $Z\in\cH$, define
\begin{equation}\label{en:eq:literal-functional}
 V(Z)=\sup_{\substack{Y\in L^2(\mathcal F;\R^d)\\\Law(Y)=\nu}}
 \E\bigl[\ip{Z+\Gamma}{Y}-\ip{Z}{X}\bigr].
\end{equation}
The first term seeks the best coupling of the Gaussian-smoothed observation
with the target law. In the second term, $Z$ and $X$ remain on the original
probability space. Their dependence must be retained: it is part of the
functional, rather than something to be chosen afresh at each evaluation.

Writing $\MCov(\rho,\nu)=\sup_{\pi\in\Pi(\rho,\nu)}\int\ip{w}{y}\,d\pi$,
attainment of the optimal transport problem gives
\begin{equation}\label{en:eq:mcov-functional}
 V(Z)=\MCov(\Law(Z)*\gamma,\nu)-\E\ip{Z}{X}.
\end{equation}
Indeed, $Z$ is $\sigma(X)$-measurable, so independence yields
$\Law(Z+\Gamma)=\Law(Z)*\gamma$. This law has an everywhere positive density
obtained by convolution with the Gaussian measure. The optimal Brenier map
$T_Z$ from this law to $\nu$ gives the admissible vector
$Y=T_Z(Z+\Gamma)$ in~\eqref{en:eq:literal-functional} and attains the maximal
covariance. The two suprema therefore agree.

Equation~\eqref{en:eq:mcov-functional} does not imply that $V(Z)$ depends only on
the law of $Z$: the term $\E\ip{Z}{X}$ still records its dependence on $X$.
In particular, we do not replace that term by $\MCov(\Law(Z),\mu)$.

The gradient is known explicitly~\citep[Lemma~2.1]{en:bps2025}:
\begin{equation}\label{en:eq:actual-gradient}
 g(Z):=\nabla V(Z)=\E[T_Z(Z+\Gamma)\mid\G]-X.
\end{equation}
This is a continuous Fr\'echet gradient on all of $\cH$, and $V$ is convex.
To see the role of the formula, fix the optimal vector $Y_Z$. It remains
admissible when $Z$ is replaced by $W$, giving
\begin{equation}\label{en:eq:global-support}
 V(W)\ge V(Z)+\ip{g(Z)}{W-Z}_{\cH}\qquad(Z,W\in\cH).
\end{equation}
Two such comparisons give the precise remainder bound
\[
0\le V(W)-V(Z)-\ip{g(Z)}{W-Z}_{\cH}
\le\norm{g(W)-g(Z)}_2\norm{W-Z}_2.
\]
Continuity of the selected optimal vector, and hence of $g$, gives the
derivative without differentiating the optimizer with respect to a parameter.
For a compact target, this continuity follows from stability of optimal plans
and uniqueness of the Brenier map.

Convex order implies equality of barycenters,
$\E X=\int y\,\nu(dy)$: apply it to both signs of each linear function. Thus
\[
 V(Z+c)=V(Z),\qquad \E g(Z)=0\qquad(c\in\R^d).
\]
Strong convexity on all of $\cH$ is therefore impossible, since constant
translations leave the functional unchanged. The quantitative curvature we
need appears after fixing the mean. For a radius $R$ and a vector $\bar z$,
write
\[
 S_{R,\bar z}=\{Z\in\cH:\norm{Z}_\infty\le R,\ \E Z=\bar z\}.
\]
This is a closed convex subset of $\cH$. To see that the condition
$\norm{Z}_\infty\le R$ is closed, take an almost surely convergent
subsequence of any $L^2$-convergent sequence.

\subsection{The main theorem}

We consider gradient-flow time $t\in[0,\infty)$ and the equation
\begin{equation}\label{en:eq:flow}
 \dot Z_t=-g(Z_t),\qquad Z_{t=0}=Z_0.
\end{equation}
This is the time parameter for motion along the functional. It is distinct
from Bass martingale time, which usually runs from $0$ to $1$.

\begin{theoremen}[Multidimensional exponential convergence]\label{en:thm:main}
Under assumptions (A1)--(A4), for every $Z_0\in\cH\cap L^\infty$ there is a
unique locally absolutely continuous solution of~\eqref{en:eq:flow} satisfying
the equation for almost every time. It is a $C^1$ curve in $\cH$, preserves
the mean, and is uniformly bounded in $L^\infty$ on the whole half-line.
There exist a global minimizer $Z^\star$ of $V$ on all of $\cH$, unique among
global minimizers with mean $\E Z_0$, and a constant $\kappa>0$ such that,
for all $t\ge0$,
\begin{align}
 \norm{Z_t-Z^\star}_2
 &\le\norm{Z_0-Z^\star}_2 e^{-\kappa t},\label{en:eq:main-norm}\\
 V(Z_t)-\inf_{W\in\cH}V(W)
 &\le\left[V(Z_0)-\inf_{W\in\cH}V(W)\right]e^{-\kappa t}.
 \label{en:eq:main-energy}
\end{align}
The constant may depend on the problem data and the initial state.
\end{theoremen}

In both estimates, the prefactor is exactly the initial error. The statement
is not restricted to sufficiently large times. The initial mean is arbitrary;
no monotonicity of the coupling of $Z_0$ with $X$, finiteness of
$\supp\mu$, or smallness of $\norm{Z_0}_\infty$ is required.

\paragraph{The known result we use.}
Applied to the same functional, Theorem~1.8 of~\citep{en:bps2025} already gives
existence and uniqueness of the flow, its strong $L^2$ limit, global
minimality, and uniqueness of the minimizer at a fixed mean. For bounded
initial data, it also gives a single essential bound for the entire
trajectory. Continuity in~\eqref{en:eq:actual-gradient} upgrades the locally
absolutely continuous curve to a $C^1$ curve; the identity $\E g=0$ preserves
the mean. Uniqueness among locally absolutely continuous solutions can also
be checked directly: monotonicity of the gradient of a convex functional
makes the squared distance between two such solutions nonincreasing almost
everywhere. It is zero when the initial states agree.
We use this qualitative part as a published theorem. The new step in the
proof below is positive curvature on $S_{R,\bar z}$ and the resulting rates.

In~\citep{en:bps2025}, the measures are centered at the outset. This does not
restrict applicability: for a common barycenter $m$, replace $X,Y$ by
$X-m,Y-m$ and translate both measures by $-m$. The assumptions are preserved,
and the lifted functional is unchanged, since
\[
 \E\bigl[\ip{Z+\Gamma}{Y-m}-\ip{Z}{X-m}\bigr]
 =\E\bigl[\ip{Z+\Gamma}{Y}-\ip{Z}{X}\bigr]-\ip{\E\Gamma}{m}
 =\E\bigl[\ip{Z+\Gamma}{Y}-\ip{Z}{X}\bigr].
\]
Moreover, $\sigma(X-m)=\sigma(X)$, so the space $\cH$ is the same.

This order of argument assumes neither a rate nor strong convexity in
advance. The known theorem first gives a radius $R$ bounding the flow and its limit.
We then prove a quantitative estimate on the corresponding bounded set. Section~\ref{en:sec:rates} will complete the proof of
Theorem~\ref{en:thm:main}.

\section{An entropic comparison and its endpoints}\label{en:sec:ent:comparison}

To obtain a convergence rate, we need a positive lower bound on the curvature
of the functional. The optimal transport map, however, is defined by a
variational problem, and differentiating it directly with respect to the
source random vector is inconvenient. We therefore begin with a regularized
problem in which the required derivative is expressed through an ordinary
variance.

Let \(D(\pi\Vert\rho\otimes\nu)\) denote the relative entropy of a coupling
\(\pi\) with respect to the product of its marginals. For \(\varepsilon>0\),
set
\[
 E_\varepsilon(\rho,\nu)=
 \sup_{\substack{\pi\in\Pi(\rho,\nu)\\
                  D(\pi\Vert\rho\otimes\nu)<\infty}}
 \left\{\int\ip{w}{y}\,\pi(dw,dy)
       -\varepsilon D(\pi\Vert\rho\otimes\nu)\right\}.
\]
Here \(\Pi(\rho,\nu)\) is the set of all couplings of \(\rho\) and \(\nu\).
Since the target vector is bounded, a finite first moment of \(\rho\) suffices
for the covariance term to be finite.

\subsection{Row and column potentials and the two marginal equations}

Let \(C=\operatorname{supp}\nu\subset\overline B_M(0)\).
A target potential provides an upper bound on the value of the problem.
We first choose it arbitrarily, and then select it so that the bound is
attained at the endpoint states we need.
For a bounded measurable function \(a\colon C\to\R\), define
\[
 v_a(w)=\varepsilon\log\int_C
 e^{(\ip{w}{y}-a(y))/\varepsilon}\,\nu(dy),
 \qquad
 p_{a,w}(dy)=
 e^{(\ip{w}{y}-a(y)-v_a(w))/\varepsilon}\,\nu(dy).
\]
We call \(a\) the column potential and \(v_a\) its associated row potential.
By definition, \(p_{a,w}\) is a probability measure.
Adding a constant \(c\) to the column subtracts that constant from the row
and leaves the coupling density unchanged. We use this freedom to normalize
the row by \(v_a(0)=0\).

For any such column and any source measure $\rho$ with a finite first moment,
the measure $Q_a(dw,dy)=\rho(dw)p_{a,w}(dy)$ is a probability measure.
Its first marginal is $\rho$; its second may still differ from $\nu$.
For a coupling $\pi\in\Pi(\rho,\nu)$ with finite relative entropy,
nonnegativity of $D(\pi\Vert Q_a)$ gives weak duality:
\begin{equation}\label{en:eq:ent:weak-duality}
 \int\ip{w}{y}\,d\pi-\varepsilon D(\pi\Vert\rho\otimes\nu)
 \le\int v_a\,d\rho+\int a\,d\nu.
\end{equation}
Indeed, the logarithm of the density of $Q_a$ with respect to
$\rho\otimes\nu$ is $(\ip{w}{y}-a(y)-v_a(w))/\varepsilon$.
Substituting it into $D(\pi\Vert Q_a)\ge0$ yields the inequality.
Every term is integrable because the target and the column are bounded,
$\rho$ has a finite first moment, and the row is Lipschitz.
This is why the upper comparison is also valid for an intermediate column
that is not optimal.

\begin{lemmaen}\label{en:lem:ent:calculus}
The function \(v_a\) is twice continuously differentiable, convex, and
\(M\)-Lipschitz. For any \(w,q,q'\in\R^d\),
\begin{align}
 Dv_a(w)[q]&=\int\ip{q}{y}\,p_{a,w}(dy),\label{en:eq:ent:first-derivative}\\
 D^2v_a(w)[q,q']&=\varepsilon^{-1}
 \operatorname{Cov}_{p_{a,w}}\bigl(\ip{q}{y},\ip{q'}{y}\bigr).
 \label{en:eq:ent:hessian}
\end{align}
\end{lemmaen}
\begin{proof}
On every bounded set of values of \(w\), the exponential in the integrand
and each of its derivatives are bounded by a common constant times a power
of \(\norm{y}\). Since \(\norm{y}\le M\) and \(\nu\) is a probability measure,
these bounds are integrable. We may therefore differentiate the integral
successively. Differentiating the logarithm gives the first formula;
differentiating once more gives the second moment minus the product of the
first moments, which is the covariance. It is nonnegative when \(q'=q\),
so the row is convex. Finally, its gradient is an average of vectors in
\(\overline B_M(0)\), and hence has norm at most \(M\).
\end{proof}

We need more than these formally defined rows: we need pairs that actually
solve the regularized problem. The next proposition supplies equality at
both endpoints of the segment that we will later interpolate.

\begin{propositionen}\label{en:prop:ent:endpoints}
Let \(Z\) be a bounded random vector independent of a standard Gaussian
vector \(\Gamma\), and let \(W=Z+\Gamma\) and \(\rho=\Law(W)\).
For every \(\varepsilon>0\), there exists a column potential \(a\), bounded
and continuous on \(C\) and convex on \(C\), such that
\begin{align}
 v_a(0)&=0,\label{en:eq:ent:normalization}\\
 \E e^{(\ip{W}{y}-v_a(W))/\varepsilon}
   &=e^{a(y)/\varepsilon}\qquad(y\in C).\label{en:eq:ent:column}
\end{align}
The measure
\[
 \pi_a(dw,dy)=
 e^{(\ip{w}{y}-a(y)-v_a(w))/\varepsilon}\,\rho(dw)\nu(dy)
\]
has marginals \(\rho,\nu\), has finite relative entropy, and attains
\(E_\varepsilon(\rho,\nu)\). In particular,
\begin{equation}\label{en:eq:ent:primal-equality}
 E_\varepsilon(\rho,\nu)=\E v_a(W)+\int a\,d\nu.
\end{equation}
\end{propositionen}

\begin{proof}
First consider a finite source
\(\rho=\sum_{i=1}^n p_i\delta_{w_i}\), where \(p_i>0\). Minimize
\[
 G(b)=\sum_i p_i b_i+
 \varepsilon\int_C\log\left(
 \sum_i p_i e^{(\ip{w_i}{y}-b_i)/\varepsilon}\right)\nu(dy)
 \quad(b\in\R^n).
\]
Adding the same constant to all \(b_i\) leaves this function unchanged.
Choose the normalization \(\min_i b_i=0\). If
\(\lvert\ip{w_i}{y}\rvert\le K\) on \(C\) and
\(p_{\min}=\min_i p_i\), then
\[
 G(b)\ge p_{\min}\max_i b_i-K+\varepsilon\log p_{\min}.
\]
Indeed, the first term is at least \(p_{\min}\max_i b_i\), and in the sum
inside the logarithm we may retain an index with \(b_i=0\). Thus the closed
sublevel sets under this normalization are compact, and the minimum is
attained. Invariance under common shifts means that the vector found is also
an unconstrained minimizer.

Set
\[
 a(y)=\varepsilon\log\sum_i p_i
 e^{(\ip{w_i}{y}-b_i)/\varepsilon}.
\]
The first-order condition for each \(b_i\) gives
\[
 \int_C e^{(\ip{w_i}{y}-a(y))/\varepsilon}\,\nu(dy)
 =e^{b_i/\varepsilon}.
\]
Hence \(v_a(w_i)=b_i\). The definition of \(a\) gives the second marginal
equation. A common shift of the pair normalizes \(v_a(0)=0\).
The function \(a\) is convex and continuous, since it is the logarithm of an
integral of exponentials of affine functions.

Now return to \(W=Z+\Gamma\). Choose finite-valued simple approximations
\(W_n\to W\) almost surely and in \(L^2\), with
\(\norm{W_n}\le2\norm{W}\). For each approximation, construct a pair
\((a_n,v_n)\) with \(v_n(0)=0\). Lemma~\ref{en:lem:ent:calculus} gives
\(\lvert v_n(w)\rvert\le M\norm{w}\).

We next show that, for fixed \(\varepsilon\), the columns are uniformly
bounded and equicontinuous. Their second marginal equation reads
\[
 a_n(y)=\varepsilon\log
 \E e^{(\ip{W_n}{y}-v_n(W_n))/\varepsilon}.
\]
On \(C\), the absolute value of the exponent is at most
\(2M\norm{W_n}/\varepsilon\le4M\norm{W}/\varepsilon\).
The norm of \(W\) has exponential moments of every order, since \(Z\) is
bounded and \(\Gamma\) is Gaussian. The integral therefore has an upper
bound uniform in \(n,y\). Jensen's inequality gives the lower bound
\[
 a_n(y)\ge-2M\E\norm{W_n}\ge-4M\E\norm{W}.
\]
As \(y\) varies, the derivative of the integral is bounded in norm by
\(\E[\norm{W_n}e^{2M\norm{W_n}/\varepsilon}]/\varepsilon\).
Together with the positive lower bound on the integral itself, this gives
a common Lipschitz bound for \(a_n\) on \(C\). All the moments used here are
finite and bounded independently of \(n\).

By the Arzel\`a--Ascoli theorem, pass to a subsequence such that
\(a_n\to a\) uniformly on \(C\). The definition of the row gives
\[
 \sup_{w\in\R^d}\lvert v_{a_n}(w)-v_a(w)\rvert
 \le\sup_{y\in C}\lvert a_n(y)-a(y)\rvert.
\]
Thus \(v_a(0)=0\), and \(v_n(W_n)\to v_a(W)\) almost surely.
The same Gaussian exponential majorant allows passage to the limit in the
column equation. We obtain~\eqref{en:eq:ent:column} for every \(y\in C\).
The limiting column is convex and continuous on \(C\).

The first marginal equation holds by the definition of the row, and the
second follows from~\eqref{en:eq:ent:column}. Hence \(\pi_a\) has the required
marginals. The logarithm of its density is integrable with respect to
\(\pi_a\): the column is bounded, \(\lvert v_a(w)\rvert\le M\norm{w}\),
and \(\lvert\ip{w}{y}\rvert\le M\norm{w}\). Its entropy is therefore finite.

It remains to verify optimality of this particular pair. For any coupling
\(\pi\) with finite relative entropy, nonnegativity of
\(D(\pi\Vert\pi_a)\) gives
\[
 \int\ip{w}{y}\,d\pi-\varepsilon D(\pi\Vert\rho\otimes\nu)
 \le \int v_a\,d\rho+\int a\,d\nu .
\]
All the affine and row terms are integrable, so subtracting the logarithm
of the density is justified. Equality holds at \(\pi=\pi_a\). This proves
both attainment and~\eqref{en:eq:ent:primal-equality}.
\end{proof}

\begin{remarken}\label{en:rem:ent:same-endpoints}
We will vary the column between two such endpoint values. An intermediate
column need not be optimal. We require only weak duality in the interior
of the segment and equality at its endpoints. There will therefore be no
need to differentiate a family of optimizers.
\end{remarken}

\subsection{The comparison function}\label{en:sec:ent:scalar}

Fix bounded $Z_0,Z_1\in\cH$ and $\varepsilon>0$. For the sources
$\Law(Z_i+\Gamma)$, choose columns $a_i$ from
Proposition~\ref{en:prop:ent:endpoints}, and set
$a_s=(1-s)a_0+sa_1$, $0\le s\le1$.

Set \(Z_s=(1-s)Z_0+sZ_1\), \(h=Z_1-Z_0\), and \(b=a_1-a_0\), and consider
\[
 J(s)=\E v_{a_s}(Z_s+\Gamma)+\int_C a_s\,d\nu .
\]
Weak duality gives
\(E_\varepsilon(\Law(Z_s+\Gamma),\nu)\le J(s)\), with equality at \(s=0,1\).
Since \(h\), \(b\), and the target vector are bounded, we may differentiate
twice with respect to \(s\) under the expectation. For example, the first
derivative of the row is bounded by
\(M\norm{h}_\infty+\norm{b}_\infty\), and the second by the square of this
constant divided by \(\varepsilon\). We obtain
\begin{equation}\label{en:eq:curv:comparison-second}
 J''(s)=\frac1\varepsilon\E
 \operatorname{Var}_{p_{a_s,Z_s+\Gamma}}
       \bigl(\ip{h}{y}-b(y)\bigr).
\end{equation}

\begin{figure}[H]
\centering
\begin{tikzpicture}[x=10cm,y=4.2cm,font=\small]
  \draw[-{Stealth[length=2mm]},draw=black!65,line width=.5pt]
    (-.025,-.68) -- (1.07,-.68) node[right] {\(s\)};
  \draw[draw=black!40,line width=.5pt] (0,0) -- (1,0)
    node[midway,above=6pt] {the common affine chord has been subtracted};
  \draw[draw=black!55,densely dotted,line width=.9pt]
    plot[domain=0:1,samples=65] (\x,{-\x*(1-\x)});
  \draw[draw=geoPthreehalf,densely dashed,line width=1.05pt]
    plot[domain=0:1,samples=65] (\x,{-1.8*\x*(1-\x)});
  \draw[draw=geoPtwo,line width=1.15pt]
    plot[domain=0:1,samples=65]
      (\x,{-1.8*\x*(1-\x)-.8*\x*\x*(1-\x)*(1-\x)});
  \fill[geoPtwo] (0,0) circle[x radius=.004,y radius=.0114];
  \fill[geoPtwo] (1,0) circle[x radius=.004,y radius=.0114];
  \node[below=4pt] at (0,-.68) {\(0\)};
  \node[below=4pt] at (1,-.68) {\(1\)};
  \node[align=center,text=black!70,font=\footnotesize] at (.5,-.18)
    {\(-\frac K2s(1-s)\)};
  \node[align=center,text=geoPthreehalf] at (.49,-.35)
    {\(J(s)-\ell_\varepsilon(s)\)};
  \node[align=center,text=geoPtwo] at (.52,-.565)
    {\(F_\varepsilon(s)-\ell_\varepsilon(s)\)};
  \draw[draw=black!40,line width=.4pt] (0,0) -- (0,-.68);
  \draw[draw=black!40,line width=.4pt] (1,0) -- (1,-.68);
\end{tikzpicture}
\caption{Schematic comparison for a fixed \(\varepsilon>0\).
Here \(F_\varepsilon(s)=E_\varepsilon(\Law(Z_s+\Gamma),\nu)\),
and \(J(s)\) is the upper comparison function built from the interpolated
target potentials. The common chord has been subtracted:
\(\ell_\varepsilon(s)=(1-s)F_\varepsilon(0)+sF_\varepsilon(1)\).
If \(J''\ge K\) for some \(K>0\), then
\(F_\varepsilon(s)-\ell_\varepsilon(s)
\le J(s)-\ell_\varepsilon(s)\le-(K/2)s(1-s)\).
All three expressions vanish at the endpoints.
The curves illustrate the directions of the inequalities; they are not
computed values of the entropic functional.}
\label{en:fig:bass-comparison}
\end{figure}
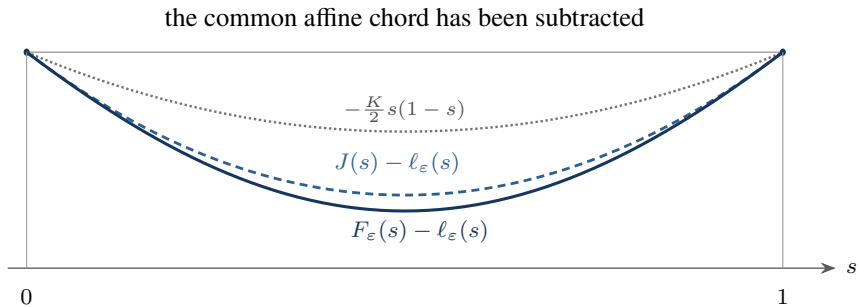

\subsection{Removing the entropy regularization}

The positive lower bound on curvature must survive as
\(\varepsilon\downarrow0\). First, we justify the passage to the
unregularized covariance problem itself.

\begin{propositionen}\label{en:prop:ent:removal}
If \(\rho\) is a probability measure with a finite first moment
and \(\nu\) has compact support, then
\[
 E_\varepsilon(\rho,\nu)\longrightarrow\MCov(\rho,\nu)
 \qquad(\varepsilon\downarrow0).
\]
\end{propositionen}
\begin{proof}
The upper bound \(E_\varepsilon\le\MCov\) follows from nonnegativity of
relative entropy. For the reverse bound, fix \(a>0\) and partition \(C\)
into finitely many measurable cells \(C_j\) of diameter at most \(a\).
Discard cells of zero mass. For an arbitrary coupling
\(\pi\in\Pi(\rho,\nu)\), define
\[
 \pi^a=\sum_{\nu(C_j)>0}
 \pi(\,\cdot\,\times C_j)\otimes\frac{\nu|_{C_j}}{\nu(C_j)}.
\]
Both marginals are preserved. If
\(r_j=d\pi(\,\cdot\,\times C_j)/d\rho\), then \(0\le r_j\le1\).
The density of \(\pi^a\) with respect to \(\rho\otimes\nu\) on \(C_j\)
is therefore \(r_j(w)/\nu(C_j)\), which is at most
\(1/\min_{\nu(C_j)>0}\nu(C_j)\). Its entropy is bounded by a finite constant
\(H_a\) independent of the original coupling.

Under resampling, the two target points remain in the same cell, so the
change in \(\ip{w}{y}\) has absolute value at most \(a\norm{w}\). Hence
\[
 \int\ip{w}{y}\,d\pi^a
 \ge\int\ip{w}{y}\,d\pi-a\int\norm{w}\,\rho(dw).
\]
Taking the supremum over \(\pi\) gives
\[
 \MCov(\rho,\nu)-a\int\norm{w}\,d\rho-\varepsilon H_a
 \le E_\varepsilon(\rho,\nu)\le\MCov(\rho,\nu).
\]
First send \(\varepsilon\) to zero at fixed \(a\), and then send \(a\) to
zero. No uniform bound on \(H_a\) as \(a\downarrow0\) is needed.
\end{proof}

We can therefore work with smooth rows at positive \(\varepsilon\), provided
we find a bound that does not degenerate as the parameter tends to zero.
The next section establishes precisely such a bound.

\section{Positive curvature on a bounded slice}\label{en:sec:curv:main}

Fix a radius \(0\le R<\infty\). We want to compare two source vectors
\(Z_0,Z_1\in\cH\) satisfying \(\norm{Z_i}_\infty\le R\).
The lower bound should depend on this common radius, but not on the vectors
chosen within it. We first derive a geometric estimate for the row potentials
and then transfer it to the random source vectors.

\subsection{An interior column bound and growth of the row potential}

Choose \(y_c\in I\) and \(r>0\) such that
\(\overline B_{2r}(y_c)\subset I\). The density of \(\nu\) on this ball
is bounded below by some \(m>0\). Write
\(G_d=\E\norm{\Gamma}\) and \(C\subset\overline B_M(0)\).

\begin{lemmaen}\label{en:lem:curv:cap-growth}
Let \(a_0,a_1\) be the endpoint columns from
Proposition~\ref{en:prop:ent:endpoints}, corresponding to \(Z_0+\Gamma,Z_1+\Gamma\),
and let \(a_s=(1-s)a_0+sa_1\), \(0\le s\le1\).
There are constants \(c>0,K<\infty\), independent of \(s\), the endpoints,
and \(0<\varepsilon\le1\), such that
\begin{equation}\label{en:eq:curv:row-growth}
 v_{a_s}(w)\ge\ip{w}{y_c}+c\norm{w}-K,\qquad v_{a_s}(0)\le0 .
\end{equation}
All these row potentials are convex and \(M\)-Lipschitz.
\end{lemmaen}
\begin{proof}
Set \(K_0=2M(R+G_d)\). The column equation and Jensen's inequality give
\[
 a_i(y)\ge\E[\ip{Z_i+\Gamma}{y}-v_{a_i}(Z_i+\Gamma)]\ge-K_0.
\]
The endpoint equality, the bound
\(E_\varepsilon(\Law(Z_i+\Gamma),\nu)\le M\E\norm{Z_i+\Gamma}\),
and \(\lvert v_{a_i}(w)\rvert\le M\norm{w}\) imply
\(\int a_i\,d\nu\le K_0\). Both estimates are preserved under interpolation
of the columns.

For \(y\in\overline B_r(y_c)\), the ball \(B_r(y)\) is contained in
\(B_{2r}(y_c)\). The function \(a_s+K_0\) is convex and nonnegative.
Jensen's inequality for the uniform measure on the ball, followed by the
lower density bound, gives
\[
 a_s(y)+K_0\le\frac1{|B_r|}\int_{B_r(y)}(a_s+K_0)\,dy
 \le\frac1{m|B_r|}\int_C(a_s+K_0)\,d\nu
 \le\frac{2K_0}{m|B_r|}.
\]
Let \(A\) denote the resulting upper bound for \(a_s(y)\).
For \(w\ne0\), the ball of radius \(r/4\) centered at
\(y_c+(r/2)w/\norm{w}\) lies in \(B_r(y_c)\), and throughout this ball
\(\ip{w}{y}\ge\ip{w}{y_c}+(r/4)\norm{w}\).
Its mass is at least \(q=m|B_{r/4}|>0\), with \(q\le1\).
Restricting the row integral to this ball yields
\[
 v_{a_s}(w)\ge\ip{w}{y_c}+(r/4)\norm{w}-A+\varepsilon\log q
 \ge\ip{w}{y_c}+(r/4)\norm{w}-A+\log q .
\]
For \(w=0\), the same bound follows directly from the integral.
We may therefore take \(c=r/4\) and \(K=A-\log q\).

Finally, the endpoint normalizations mean that
\(\int e^{-a_i/\varepsilon}\,d\nu=1\). Convexity of the exponential gives
\[
 \int e^{-a_s/\varepsilon}\,d\nu
 \le(1-s)\int e^{-a_0/\varepsilon}\,d\nu
       +s\int e^{-a_1/\varepsilon}\,d\nu=1.
\]
Hence \(v_{a_s}(0)\le0\). Convexity and the Lipschitz bound follow from
Lemma~\ref{en:lem:ent:calculus}.
\end{proof}

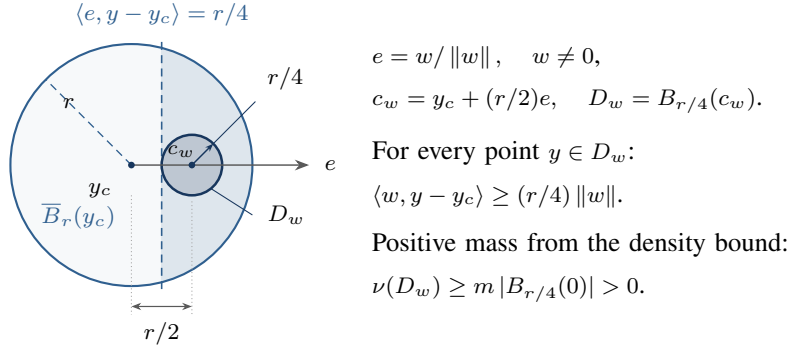
\begin{figure}[H]
\centering
\begin{tikzpicture}[font=\small,
  dim/.style={{Stealth[length=1.3mm]}-{Stealth[length=1.3mm]},draw=black!65,line width=.45pt}]
  \fill[geoPthreehalf!4] (0,0) circle[radius=1.6];
  \begin{scope}
    \clip (0,0) circle[radius=1.6];
    \fill[geoPthreehalf!16] (.4,-1.65) rectangle (1.65,1.65);
  \end{scope}
  \draw[draw=geoPthreehalf,line width=.8pt] (0,0) circle[radius=1.6];
  \draw[draw=geoPthreehalf,densely dashed,line width=.6pt] (.4,-1.64) -- (.4,1.64);
  \node[anchor=south,text=geoPthreehalf] at (.4,1.71)
    {\(\ip{e}{y-y_c}=r/4\)};
  \filldraw[draw=geoPtwo,fill=geoPtwo!28,line width=.95pt]
    (.8,0) circle[radius=.4];
  \draw[-{Stealth[length=1.9mm]},draw=black!65,line width=.55pt]
    (0,0) -- (2.35,0) node[right=2pt] {\(e\)};
  \draw[draw=geoPthreehalf,densely dashed,line width=.5pt]
    (0,0) -- (-1.13137,1.13137) node[midway,above left=2pt] {\(r\)};
  \fill[geoPtwo] (0,0) circle[radius=1.25pt];
  \node[anchor=north east] at (-.12,-.10) {\(y_c\)};
  \fill[geoPtwo] (.8,0) circle[radius=1.25pt];
  \node[font=\footnotesize] at (.66,.18) {\(c_w\)};
  \node[text=geoPthreehalf] at (-.70,-.67) {\(\overline B_r(y_c)\)};
  \draw[-{Stealth[length=1.1mm]},draw=geoPtwo,line width=.55pt]
    (.8,0) -- (1.08284,.28284);
  \draw[draw=geoPtwo,line width=.4pt] (1.08284,.28284) -- (1.64,.84)
    node[above right=1pt] {\(r/4\)};
  \draw[draw=geoPtwo,line width=.4pt] (1.05,-.31) -- (1.60,-.66)
    node[right=2pt] {\(D_w\)};
  \draw[draw=black!35,densely dotted,line width=.4pt] (0,-.40) -- (0,-1.98);
  \draw[draw=black!35,densely dotted,line width=.4pt] (.8,-.48) -- (.8,-1.98);
  \draw[dim] (0,-1.87) -- (.8,-1.87);
  \node[below=2pt] at (.4,-1.87) {\(r/2\)};
  \node[anchor=west,align=left,text width=6.7cm] at (3.05,1.13)
    {\(e=w/\norm{w},\quad w\ne0\),\\[4pt]
     \(c_w=y_c+(r/2)e,\quad D_w=B_{r/4}(c_w)\).};
  \node[anchor=west,align=left,text width=6.7cm] at (3.05,-.12)
    {For every point \(y\in D_w\):\\[4pt]
     \(\ip{w}{y-y_c}\ge(r/4)\norm{w}\).};
  \node[anchor=west,align=left,text width=6.7cm] at (3.05,-1.35)
    {Positive mass from the density bound:\\[4pt]
     \(\nu(D_w)\ge m\,|B_{r/4}(0)|>0\).};
\end{tikzpicture}
\caption{Geometry of the row-growth bound. The shaded region is the part of
\(\overline B_r(y_c)\) on the \(\ip{e}{y-y_c}\ge r/4\) side of the dashed
line. The whole small ball lies on that side, and
\(\norm{c_w-y_c}+r/4=3r/4<r\) ensures \(D_w\subset B_r(y_c)\).
The picture is two-dimensional; the argument holds in \(\R^d\).
The ambient choice \(\overline B_{2r}(y_c)\subset I\) is not drawn.}
\label{en:fig:bass-caps}
\end{figure}

We now need to extract a positive amount of curvature from the growth of
the row potential. Its Hessian need not be positive definite at every point:
an integrated lower bound on one common ball will suffice.

\begin{lemmaen}\label{en:lem:curv:hessian-mass}
Let \(v\colon\R^d\to\R\) be convex, \(C^2\), and
\(M\)-Lipschitz. Suppose that \(v(0)\le0\) and
\[
 v(w)\ge\ip{w}{y_c}+c\norm{w}-K,\qquad c>0.
\]
Then there exist \(T<\infty\) and \(\eta>0\), depending only on
\(M,y_c,c,K,d\), such that
\begin{equation}\label{en:eq:curv:hessian-mass}
 \int_{B_T(0)}D^2v(w)[q,q]\,dw\ge\eta\norm{q}^2
 \qquad(q\in\R^d).
\end{equation}
\end{lemmaen}
\begin{proof}
Fix \(\rho>0\). For \(\norm{z}<\rho\) and \(\norm{\xi}=1\),
the two lower bounds and \(v(z)\le M\norm{z}\) imply
\[
 v(z+L\xi)+v(z-L\xi)-2v(z)
 \ge2cL-2K-2(\norm{y_c}+M)\rho .
\]
Choose \(L>0\) large enough that the right-hand side is at least
\(cL/2\) for all such \(z,\xi\).
For \(j(t)=v(z+t\xi)\), convexity and the fundamental theorem of calculus give
\[
 j(L)+j(-L)-2j(0)
 \le L[j'(L)-j'(-L)]
 =L\int_{-L}^{L}D^2v(z+t\xi)[\xi,\xi]\,dt.
\]
Thus the integral on the right, without the factor \(L\), is at least
\(c/2\). Integrate this inequality over \(z\in B_\rho(0)\).
The integrand is nonnegative. For each fixed \(t\), its integral over
the translated ball is no greater than its integral over \(B_T(0)\),
provided that \(T=L+\rho+1\). Fubini's theorem gives
\[
 \frac c2|B_\rho|\le
 2L\int_{B_T(0)}D^2v(w)[\xi,\xi]\,dw.
\]
We may take \(\eta=c|B_\rho|/(4L)>0\).
For an arbitrary nonzero \(q\), apply the bound to
\(\xi=q/\norm{q}\) and use bilinearity of the Hessian.
For \(q=0\), the assertion is immediate.
\end{proof}

By Lemma~\ref{en:lem:curv:cap-growth}, the same \(T,\eta\) work for all
interpolated row potentials with \(0<\varepsilon\le1\).
This independence of \(\varepsilon\) is what will allow us to return
to the original problem.

\subsection{Gaussian observations and two copies}

Return to the comparison function $J$ from Section~\ref{en:sec:ent:scalar}.
Here $Z_s=(1-s)Z_0+sZ_1$, $h=Z_1-Z_0$, and $b=a_1-a_0$.

How can we eliminate \(b\)? Take two independent values of \(h\) conditional
on the same observation. In the difference of the corresponding expressions,
the column increment cancels.

\begin{lemmaen}\label{en:lem:curv:posterior}
For \(T,\eta\) from Lemma~\ref{en:lem:curv:hessian-mass}, set
\[
 g_T=(2\pi)^{-d/2}e^{-(R+T)^2/2},
 \qquad b_\gamma=(2\pi)^{-d/2}.
\]
Then, for every \(0\le s\le1\),
\begin{equation}\label{en:eq:curv:comparison-lower}
 J''(s)\ge\delta_R\,\E\norm{h-\E h}^2,
 \qquad \delta_R=\frac{g_T^2\eta}{2b_\gamma}>0 .
\end{equation}
\end{lemmaen}
\begin{proof}
Write \(\lambda_s=\Law(Z_s,h)\) for the joint law on
\(\R^d\times\R^d\), and let \(\varphi\) be the standard Gaussian
density. Since both \(Z_0,Z_1\) are \(\G\)-measurable,
the pair \((Z_s,h)\) is independent of \(\Gamma\).
The observation \(w=Z_s+\Gamma\) has density
\[
 \beta_s(w)=\int\varphi(w-z)\,\lambda_s(dz,dh)>0.
\]
The conditional law of the hidden pair is the explicit probability measure
\[
 \lambda_{s,w}(dz,dh)=
 \frac{\varphi(w-z)}{\beta_s(w)}\,\lambda_s(dz,dh).
\]
This formula is verified by integrating measurable test functions.
The conditional law is constructed on a finite-dimensional space;
no regular conditional distribution on the original \(\Omega\)
is needed.

Fix \(w\) and put \(p=p_{a_s,w}\).
In the Hilbert space of centered functions in \(L^2(p)\),
let \(U_h\) be the centered version of
\(\ip{h}{y}-b(y)\). For two independent conditional copies \(h,h'\),
\[
 \norm{U_h-U_{h'}}_{L^2(p)}^2
 \le2\norm{U_h}_{L^2(p)}^2+2\norm{U_{h'}}_{L^2(p)}^2 .
\]
Averaging gives
\[
 \E_{\lambda_{s,w}}\operatorname{Var}_p(\ip{h}{y}-b(y))
 \ge\frac14\E_{\lambda_{s,w}^{\otimes2}}
       \operatorname{Var}_p\ip{h-h'}{y}.
\]
The column \(b\) has cancelled. By~\eqref{en:eq:ent:hessian}, the last variance,
divided by \(\varepsilon\), equals
\(D^2v_{a_s}(w)[h-h',h-h']\).

On \(B_T(0)\), since \(\norm{z}\le R\), we have
\[
 g_T\le\varphi(w-z)\le b_\gamma,\qquad
 g_T\le\beta_s(w)\le b_\gamma.
\]
Multiplying the conditional two-copy integral by the observation density
leaves the factor
\[
 \frac{\varphi(w-z)\varphi(w-z')}{\beta_s(w)}
 \ge\frac{g_T^2}{b_\gamma}.
\]
Restrict the nonnegative integral in
\eqref{en:eq:curv:comparison-second} to \(B_T(0)\),
apply this bound, and interchange the integrals:
\[
 J''(s)\ge\frac{g_T^2}{4b_\gamma}
 \int\!\!\int\!\!\int_{B_T(0)}
 D^2v_{a_s}(w)[h-h',h-h']\,dw\,
 \lambda_s(dz,dh)\lambda_s(dz',dh').
\]
Nonnegativity allows us to use Tonelli's theorem. Finiteness of the required
integrals also follows from continuity of the Hessian on the closed ball
and boundedness of \(h\).
By~\eqref{en:eq:curv:hessian-mass}, the right-hand side is at least
\[
 \frac{g_T^2\eta}{4b_\gamma}\E\norm{h-h'}^2
 =\frac{g_T^2\eta}{2b_\gamma}\E\norm{h-\E h}^2 .
\]
The last equality is the usual variance identity for two independent
copies of a vector.
\end{proof}

\begin{figure}[H]
\centering
\begin{tikzpicture}[
  font=\small,
  bass box/.style={draw=black!40,rounded corners=2pt,align=center,
    inner sep=7pt,minimum height=1cm},
  bass arrow/.style={-{Stealth[length=2mm]},draw=black!60,line width=.65pt}
]
\node[bass box,text width=10.9cm,fill=black!3] (obs) at (0,2.5)
  {Fix \(w\): two independent draws from \(\lambda_{s,w}\).\\
   The same target law \(p_{a_s,w}\) and function \(b=a_1-a_0\)};
\node[bass box,text width=4.8cm,fill=geoPtwo!6] (left) at (-3.05,.75)
  {First conditional copy \((z,h)\)\\[3pt]
   \(A_h(y)=\ip{h}{y}-b(y)\)};
\node[bass box,text width=4.8cm,fill=geoPthreehalf!6] (right) at (3.05,.75)
  {Second conditional copy \((z',h')\)\\[3pt]
   \(A_{h'}(y)=\ip{h'}{y}-b(y)\)};
\node[bass box,text width=10.9cm,fill=geoPtwo!9] (diff) at (0,-1.1)
  {The same column increment cancels: \(b(y)-b(y)=0\).\\[4pt]
   \(A_h(y)-A_{h'}(y)=\ip{h-h'}{y}\)};
\node[align=center,text width=11.3cm] (hess) at (0,-2.65)
  {\(\varepsilon^{-1}\operatorname{Var}_{p_{a_s,w}}
       \ip{h-h'}{y}
     =D^2v_{a_s}(w)[h-h',h-h']\)};
\draw[bass arrow] (obs.south) -- ++(0,-.25) -| (left.north);
\draw[bass arrow] (obs.south) -- ++(0,-.25) -| (right.north);
\draw[bass arrow] (left.south) -- ++(0,-.3) -| ([xshift=-2.4cm]diff.north);
\draw[bass arrow] (right.south) -- ++(0,-.3) -| ([xshift=2.4cm]diff.north);
\draw[bass arrow] (diff.south) -- (hess.north);
\end{tikzpicture}
\caption{What the two-copy argument removes. For a fixed observation \(w\),
take two independent copies of the Euclidean pair from the same
conditional law \(\lambda_{s,w}\).
The difference of target potentials is the same in both expressions
and cancels. The variance of the remaining linear expression is given
by the Hessian of the log-integral.}
\label{en:fig:bass-two-copies}
\end{figure}
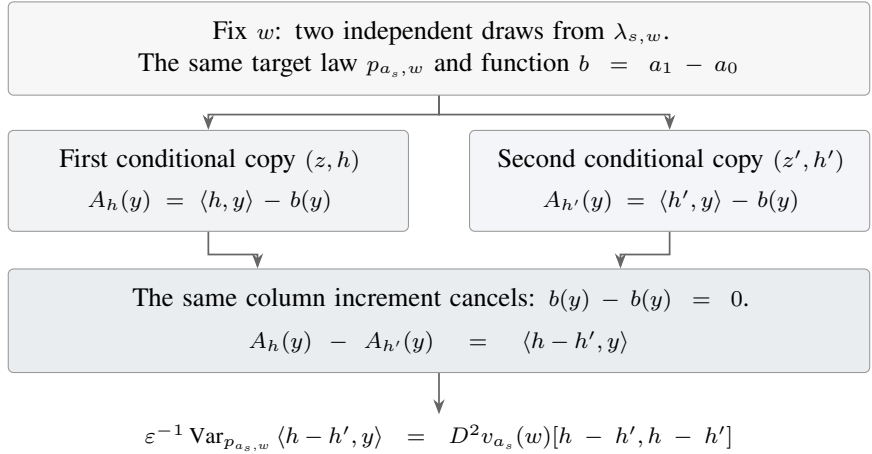

\subsection{Returning to the Bass functional}

\begin{theoremen}\label{en:thm:curv:bounded-slice}
For every \(0\le R<\infty\), there exists \(m_R>0\) such that whenever
\(Z,W\in\cH\), \(\norm{Z}_\infty,\norm{W}_\infty\le R\), and
\(\E Z=\E W\),
\begin{equation}\label{en:eq:curv:strong-support}
 V(W)\ge V(Z)+\ip{g(Z)}{W-Z}_{\cH}
             +\frac{m_R}{2}\norm{W-Z}_2^2 .
\end{equation}
The constant is chosen independently of \(Z,W\) and their common mean.
\end{theoremen}
\begin{proof}
For the two endpoints \(Z_0,Z_1\), Lemma~\ref{en:lem:curv:posterior} gives
\(J''(s)\ge\delta_R\E\norm{h-\E h}^2\).
Consequently,
\[
 J(s)\le(1-s)J(0)+sJ(1)
 -\frac{\delta_R}{2}s(1-s)\E\norm{h-\E h}^2 .
\]
This is ordinary convexity of the function
\(J(s)-(\delta_R/2)s^2\E\norm{h-\E h}^2\).
Use \(E_\varepsilon(\Law(Z_s+\Gamma),\nu)\le J(s)\)
on the left-hand side, and the endpoint equalities
from~\eqref{en:eq:ent:primal-equality} on the right-hand side. Subtracting
\(\E\ip{Z_s}{X}\) leaves the curvature unchanged, since this term
is affine in \(s\).

The constant \(\delta_R\) is independent of \(0<\varepsilon\le1\).
By Proposition~\ref{en:prop:ent:removal}, letting
\(\varepsilon\downarrow0\) gives
\[
 V(Z_s)\le(1-s)V(Z_0)+sV(Z_1)
 -\frac{\delta_R}{2}s(1-s)\E\norm{h-\E h}^2 .
\]
If the endpoint means agree, then \(\E h=0\), and the last expectation
equals \(\norm{Z_1-Z_0}_2^2\).
Finally, set \(Z_0=Z\), \(Z_1=W\), subtract \(V(Z)\),
divide by \(s>0\), and let \(s\downarrow0\).
The Fr\'echet derivative of the functional yields
\eqref{en:eq:curv:strong-support} with \(m_R=\delta_R\).
\end{proof}

\noindent\textit{Comment and intuition.}
The functional cannot be strongly convex on the whole space:
a common translation of the source vector preserves its value because
\(\mu\) and \(\nu\) have the same barycenter. The centered increment
\(h-\E h\) therefore appears naturally in the estimate.
Fixing the mean removes precisely this degeneracy.
Gaussian noise then allows us to compare any two bounded configurations:
on one common observation ball, their densities have a common positive
lower bound. The geometry of the target measure supplies a positive
amount of curvature on this ball. Together, these two facts give
the positive constant \(m_R\).

\section{From curvature to convergence rates}\label{en:sec:rates}

We have obtained a quantitative bound for two bounded states with the same
mean. Why can it now be applied at all times?
The qualitative theorem recalled in Section~\ref{en:sec:setting} supplies one
$R<\infty$ such that $\norm{Z_t}_\infty\le R$ for every $t\ge0$.
The strong $L^2$ limit $Z^\star$ belongs to the same closed set.
All these vectors have mean $\bar z=\E Z_0$.
Thus one fixed strong convexity constant $m_R>0$ applies throughout
the argument below; it does not have to be chosen anew at each time.

Recall the bound:
\begin{equation}\label{en:eq:rate-support}
 V(W)\ge V(Z)+\ip{g(Z)}{W-Z}_{\cH}
       +\frac{m_R}{2}\norm{W-Z}_2^2\qquad(Z,W\in S_{R,\bar z}).
\end{equation}
The gradient vanishes at $Z^\star$: $g(Z^\star)=0$.
This follows from differentiability of $V$ and minimality of $Z^\star$
on the whole space $\cH$.

\subsection{Distance to the limit}

Interchange $Z,W$ in~\eqref{en:eq:rate-support} and add the two inequalities.
The functional values cancel, leaving strong monotonicity:
\[
 \ip{g(Z)-g(W)}{Z-W}_{\cH}\ge m_R\norm{Z-W}_2^2.
\]
Taking $Z=Z_t$ and $W=Z^\star$ gives
\begin{align*}
 \frac{d}{dt}\norm{Z_t-Z^\star}_2^2
 &=-2\ip{g(Z_t)-g(Z^\star)}{Z_t-Z^\star}_{\cH}\\
 &\le-2m_R\norm{Z_t-Z^\star}_2^2.
\end{align*}
The derivative exists because the flow is a $C^1$ curve.
Multiplication by $e^{2m_Rt}$ shows that
$e^{2m_Rt}\norm{Z_t-Z^\star}_2^2$ is nonincreasing. At $t=0$,
its value is the squared initial distance. Taking square roots yields
\begin{equation}\label{en:eq:strong-norm-rate}
 \norm{Z_t-Z^\star}_2\le e^{-m_Rt}\norm{Z_0-Z^\star}_2.
\end{equation}
This argument also covers $Z_0=Z^\star$: both sides then vanish.

The estimate also controls the laws in the quadratic Wasserstein distance.
The joint law of $(Z_t,Z^\star)$ is an admissible coupling of their
respective laws, so the definition of that distance gives
\[
 W_2(\Law(Z_t),\Law(Z^\star))\le\norm{Z_t-Z^\star}_2.
\]
Control of the original Hilbert-space trajectory therefore also gives
a quantitative approximation of the desired distribution.

\subsection{The functional value gap}

Decay of the distance alone does not yet provide the required value estimate
with the exact initial prefactor. We therefore consider separately
$F(t)=V(Z_t)-V(Z^\star)\ge0$.
Take $Z=Z_t$, $W=Z^\star$ in~\eqref{en:eq:rate-support}. Writing
$r_t=\norm{Z_t-Z^\star}_2$, Cauchy--Schwarz gives
\[
 F(t)\le\norm{g(Z_t)}_2r_t-\frac{m_R}{2}r_t^2
 \le\frac{\norm{g(Z_t)}_2^2}{2m_R}.
\]
The last step is completion of the square:
$ar-m_Rr^2/2=a^2/(2m_R)-m_R(r-a/m_R)^2/2$.
No minimization outside the slice is used: the supporting inequality
was applied to the specific admissible vector $Z^\star$.

The chain rule and the flow equation give
\[
 F'(t)=\ip{g(Z_t)}{\dot Z_t}_{\cH}
       =-\norm{g(Z_t)}_2^2\le-2m_RF(t).
\]
The same integrating-factor argument yields the stronger estimate
\begin{equation}\label{en:eq:strong-energy-rate}
 V(Z_t)-V(Z^\star)
 \le e^{-2m_Rt}\bigl[V(Z_0)-V(Z^\star)\bigr].
\end{equation}
Since $V(Z^\star)=\inf_{\cH}V$, choosing $\kappa=m_R$ gives
both~\eqref{en:eq:main-norm} and~\eqref{en:eq:main-energy} simultaneously.
All other conclusions of the main theorem are already supplied by
the qualitative result in Section~\ref{en:sec:setting}.
This completes its proof.

\paragraph{Comment and intuition.}
The qualitative theory keeps the trajectory in one region. On that region,
uniform curvature prevents two states with the same mean from remaining
far apart while their gradients differ only slightly.
The distance to the minimizer therefore decays with a single exponential rate.
For the energy, we also used the relation between its gap and the squared
norm of the gradient. This gives two different differential inequalities
and two exact initial prefactors.

We do not claim that $m_R$ or the rate exponent is optimal.
The constants in the preceding section were chosen to establish positivity
and uniformity; their numerical values may be far from the best ones.

\section{Significance and relation to other work}
\label{en:sec:related}

\subsection{The question answered by the estimate}

In martingale optimal transport, two laws must be coupled while preserving
the conditional mean. The variational approach to Bass martingales leads
to the minimization of a lifted functional. A natural way to study this
problem is to follow its continuous gradient descent. The existence of a
limit and the rate of approach to it answer different questions: the former
identifies where the flow goes, while the latter quantifies the remaining
error in terms of the elapsed time.

Backhoff-Veraguas, Pammer, and Schachermayer~\citep{en:bps2025} established
the qualitative result in arbitrary dimension and the two exponential
estimates in dimension one. They formulated the extension of the latter
statement to higher dimensions as a conjecture. Theorem~\ref{en:thm:main}
gives this extension under the same assumptions on the data. We explicitly
use the known qualitative part; the quantitative contribution is
multidimensional strong convexity on bounded slices with a fixed mean.

\subsection{The continuous flow and discrete methods}

Hasenbichler, Joseph, Loeper, Obłój, and Pammer~\citep{en:martingaleSinkhorn}
study the martingale Sinkhorn algorithm. Its iterations construct a Bass
potential in arbitrary dimension. This is a closely related calibration
problem, but the sequence of nonlinear updates differs from the
flow~\eqref{en:eq:flow}. A convergence theorem for these iterations does not
automatically determine the rate of the solution to our differential
equation.

Marini~\citep{en:marini2026} considers the broader class of $q_x$-Bass
martingales and the associated variational and iterative constructions.
That work enlarges the family of martingale couplings under study.
Here the Gaussian smoothing, the functional~\eqref{en:eq:literal-functional},
and its continuous time parameter are fixed. We do not extend our estimates
to the entire $q_x$ class without a separate proof.

\subsection{Positive entropy and the vanishing-entropy limit}

The connection with entropic problems is particularly relevant here,
since entropy enters the proof itself. Hasenbichler, Pammer, and
Thonhauser~\citep{en:hasenbichler2026} study the Schrödinger--Bass bridge,
its variational description, and its limiting regimes. These results
explain why regularized problems arise naturally alongside the Bass
problem. Convergence of their values or optimizers alone, however, does
not provide a rate for a particular unregularized trajectory.

Quantitative analysis of entropic martingale algorithms is developed,
in particular, by Kazeykina, Ren, and Wang~\citep{en:kazeykina2026} and by
Bifronte, Ley, and Sylvestre~\citep{en:bifronte2026}. The former establish
exponential convergence of the martingale Sinkhorn algorithm in relative
entropy. The latter prove linear convergence of SISTA in arbitrary
dimension by combining control of the iterates, convexity modulo affine
symmetries, and a Polyak--Łojasiewicz-type estimate. These works help
explain the general mechanism of boundedness combined with quantitative
curvature.

The distinguishing step in our argument is the passage
$\varepsilon\downarrow0$. The second derivative of the comparison function
has a single positive lower bound for $0<\varepsilon\le1$. Strong
convexity therefore passes to the original functional. A rate for a problem
with a fixed positive parameter alone would not suffice: its constant
could vanish in the limit. Here uniformity is proved using the Gaussian
conditional distribution and the integrated Hessian mass on a single
fixed ball.

\subsection{Why integrated curvature is enough}

In multidimensional transport, directly analyzing changes in the Brenier
map can lead to regularity questions for the nonlinear transport problem.
Our proof uses a different object: a logarithmic integral whose Hessian is
expressed as a covariance. The two-copy argument eliminates the increment
of the column potential. An integrated lower bound for the Hessian then
suffices; a positive pointwise bound throughout the whole space is not
needed. Consequently, the proof requires no additional smoothness of the
density or target boundary. This describes the method used here; we make
no claim that other analytic approaches are impossible.

\paragraph{Limitations.}
The main theorem concerns bounded initial states and a compact convex
target whose density is bounded above and bounded away from zero almost
everywhere on the interior of its support, under the stated irreducibility
assumption. The qualitative theory allows more general initial states,
but we have not proved the exponential estimates obtained here for all
such states. We also do not establish a uniform rate over all pairs of
measures, an optimal exponent, a discretization rate, or the computational
complexity of calibration. This literature discussion distinguishes the
respective results without claiming to be exhaustive.

\section{Formal verification and reproducibility}\label{en:sec:lean}

The final multidimensional statement has been formalized in Lean~4 using
mathlib. The verification concerns the precise original data, functional,
and differential equation, rather than an abstract implication from strong
convexity to convergence. In particular, the formal chain derives
quantitative convexity, coercivity, a radial estimate, and boundedness of
the flow from the original assumptions. Neither an already constructed
minimizer nor an a priori bounded trajectory is supplied as an additional
hypothesis to the final theorem.

The main formal result is named
\begin{center}\small
\texttt{BassFlow.EntropicRoute.original\_Bass\_G}.
\end{center}
The entry point \texttt{actual\_sourceBass\_G} expresses the same result
in terms of the measure data. It retains the original probability space,
the law of $X$, independent Gaussian noise, and an arbitrary bounded
$\sigma(X)$-measurable initial state. Its conclusion includes both estimates
with the exact initial prefactors, minimization over the entire original
$L^2$ space, and uniqueness of the global minimizer with a prescribed mean.

The complete verifiable chain contains 233 project modules. The sources
and instructions are available in the standalone repository
\url{https://github.com/Pure-Maths/Bass-Flow-Verification}.
The proof version corresponding to this edition of the paper is pinned
to the following commit:
\begin{center}\small
\href{https://github.com/Pure-Maths/Bass-Flow-Verification/tree/c8c9dec30e36585963df3f4a69508d30c26ff51f}{\texttt{c8c9dec30e36585963df3f4a69508d30c26ff51f}}.
\end{center}
A snapshot of this pinned version is also supplied as ancillary code.
The nine external modules are reconstructed by the bootstrap step below.

The files in \texttt{Numbered} expose the checked formal statements and
their axiomatic dependencies for comparison with the numbering of the
paper. Dependency versions are fixed in \texttt{lean-toolchain} and
\texttt{lake-manifest.json}.

A fresh copy is prepared and checked with the commands
\begin{quote}\small\ttfamily
python3 scripts/bootstrap\_external.py\\
lake exe cache get\\
python3 scripts/check\_sources.py\\
python3 reproduce.py
\end{quote}
The first step reconstructs nine external modules from a pinned source
and checks their hashes before and after the compatible API port. The
other 224 project modules are included in the repository. The main entry
points relate to the mathematical proof as follows:

\begin{center}
\small
\renewcommand{\arraystretch}{1.2}
\begin{tabularx}{\linewidth}{@{}>{\raggedright\arraybackslash}p{.36\linewidth}>{\raggedright\arraybackslash}X@{}}
\toprule
Mathematical step & Main formalization modules\\
\midrule
Fréchet gradient of the original $V$ &
\texttt{ActualSelectedContinuity}\\
Entropic endpoint optimizers and the limiting inequality &
\texttt{ActualEntropyChord} and its dependencies\\
Strong convexity on a bounded slice &
\texttt{ActualBassStrongConvexity}\\
Independent derivation of boundedness and a stationary limit &
\texttt{SourceCenteredL1Coercivity},
\texttt{ActualBoundedSourceRadial},
\texttt{ActualBassStationaryFromRadial}\\
Minimality, uniqueness, and the complete result &
\texttt{ActualBassAllSourceUniqueness}, \texttt{ActualBassG}\\
\bottomrule
\end{tabularx}
\end{center}

In the paper, the qualitative part is explicitly taken from the published
Theorem~1.8 of~\citep{en:bps2025}. The formal project follows a longer,
independent chain: it also checks the construction of the flow, its
confinement to a bounded set, and minimality of the limit over all of
$\cH$. This difference in the order of presentation does not change the
final statement. Comparison with the code requires examining the formal
types and dependencies; the table is a guide, not a substitute for that
comparison.

The original package underwent two separate audits of the verifiable
chain. Their reports list only the standard Lean axioms
\texttt{propext}, \texttt{Classical.choice}, and \texttt{Quot.sound} for
the final declarations. Hash checks, recompilation of the project
modules, and trust in the pinned external library are distinct levels of
reproducibility. Commands and details of the individual runs are provided
in the README of the standalone project.

Formal verification does not establish novelty or replace scientific
peer review. Nor does it automatically formally verify every explanatory
sentence or diagram in the paper.

\section{Use of AI in preparing the paper}\label{en:sec:ai}

Language models were used in developing the proof and its formal
representation. They assisted with technical steps, comparison of possible
arguments, and preparation of Lean code. In preparing the Russian
manuscript, AI was also used to structure the text, find and compare
primary sources, prepare the LaTeX sources and illustrations, and carry
out independent checks of individual parts of the proof.

AI was also used for the English translation and its comparison with the Russian text.

Several agents with different tasks were involved in the editorial
review: mathematical checking, comparison with the formal statements,
expository order, and technical reproducibility. These checks help identify
errors and ambiguities, but do not replace scientific peer review. Lean
kernel checking concerns precise formal statements and their dependencies;
their correspondence with the ordinary mathematical text requires a
separate comparison.

AI is not listed as an author. The choice of final formulations, assessment
of the mathematical content, and responsibility for the published version
rest with the author.

\section{Conclusion}\label{en:sec:conclusion}

We have studied the convergence rate of the continuous gradient flow of
the lifted Bass functional. Under the original irreducibility assumptions,
a compact convex target support, a density bounded above and bounded away
from zero almost everywhere on the interior of the support, and a bounded
initial state, we obtained multidimensional exponential estimates for the
$L^2$ distance and the functional gap. A single positive constant applies
at every nonnegative time along the given trajectory. The prefactors in
the estimates are exactly the corresponding initial errors.

The main quantitative step is strong convexity on sets with a common
$L^\infty$ bound and a fixed mean. The entropic comparison expresses curvature as a variance.
Comparing two independent copies eliminates the change in the target
potential, while an integrated Hessian estimate yields a positive lower
bound that survives as the entropic regularization is removed. The known qualitative
theorem places the entire trajectory in one bounded set. The same
curvature constant can therefore be used over the whole time interval,
giving the required rate.

Minimization is over the entire original space of random vectors, whereas
uniqueness concerns global minimizers with a prescribed mean. These two
qualifications have different roles: the first specifies the comparison
class, and the second accounts for the symmetry of the functional.

We also record the limitations of the result. We have not identified the
optimal rate exponent or proved the existence of one rate valid for all
admissible measures and initial states. Boundedness of the initial state
and the geometric assumptions on the target measure remain hypotheses of
the theorem. Smoothness of the density or of the boundary of its support,
however, is not required.

\clearpage
\setcounter{section}{0}
\setcounter{equation}{0}
\setcounter{figure}{0}
\setcounter{table}{0}
\setcounter{footnote}{0}
\selectlanguage{russian}
\color{ink}
\makerussiantitle
\vspace{-1.2em}

\begin{abstract}
Градиентный поток может сходиться к минимуму, не имея заранее известной
скорости сходимости. Мы изучаем этот вопрос для поднятого функционала Басса
в мартингальном оптимальном транспорте. Пусть исходная и целевая меры находятся
в выпуклом порядке и образуют неприводимую пару. Предположим, что поддержка
целевой меры компактна и выпукла, а её плотность почти всюду на внутренности
поддержки ограничена сверху и отделена от нуля. Пусть поддержка исходной меры
лежит строго внутри целевой. При этих условиях доказываем многомерные
экспоненциальные оценки для ограниченного начального состояния.
Одна положительная константа даёт при всех неотрицательных временах как
убывание расстояния в $L^2$ до глобального минимизатора, так и убывание зазора
функционала. В обеих оценках множитель равен соответствующей начальной ошибке.

Основной шаг --- сильная выпуклость на множествах с единой верхней оценкой
$L^\infty$-нормы и фиксированным средним. Чтобы получить её, вводим энтропийную задачу и интерполируем
оптимальные потенциалы двух концевых состояний. Вторая производная
сравнивающей функции выражается через дисперсию. Приём с двумя копиями
устраняет неизвестное изменение потенциала, а внутренняя геометрия цели
и гауссовское усреднение дают положительную оценку, равномерную по
энтропийному параметру. После удаления энтропии она передаётся исходному
функционалу. Известная качественная теорема помещает всю траекторию в одно
ограниченное множество, что позволяет получить обе скорости.

Результат продолжает на размерности $d\ge2$ экспоненциальную теорему
Backhoff-Veraguas, Pammer и Schachermayer. Гладкости плотности или границы цели
не требуется. Итоговый результат сопровождается формальной верификацией
в Lean~4; использование ИИ описано в разделе~\ref{ru:sec:ai}.
\end{abstract}

\medskip
\noindent\textbf{Ключевые слова.} Функционал Басса, мартингальный оптимальный
транспорт, градиентный поток, энтропия, сильная выпуклость, экспоненциальная сходимость.

\noindent\textbf{Математическая классификация 2020.}
60G44, 49Q22, 34G20.

\section{Введение}\label{ru:sec:intro}

Градиентный поток уменьшает значение функционала. Но само уменьшение ещё
не говорит, насколько быстро траектория приближается к минимизатору.
Даже после доказательства сходимости остаётся отдельный вопрос:
можно ли оценить ошибку через время и начальное данное?

В этой работе мы рассматриваем этот вопрос для поднятого функционала
Басса, возникающего в мартингальном оптимальном транспорте.
В работе Backhoff-Veraguas, Pammer и Schachermayer~\citep{ru:bps2025}
для этого потока установлена сильная сходимость в \(L^2\) в нескольких
измерениях при соответствующих условиях на исходные меры.
Экспоненциальная скорость получена там в одномерном случае.
Мы доказываем многомерные экспоненциальные оценки при исходных
геометрических и вероятностных предпосылках.

Зачем такая задача возникает в мартингальном транспорте? Заданы начальное
и конечное распределения случайного процесса. Мартингальное условие требует,
чтобы условное среднее будущего состояния совпадало с текущим состоянием.
Даже при фиксированных распределениях допустимых связок может быть много;
поэтому требуется принцип, выделяющий конкретную конструкцию.
Для мартингала Басса поиск подходящего исходного закона сводится к
вариационной задаче, а поднятый функционал позволяет изучать её в
гильбертовом пространстве~\citep{ru:bps2025}. Именно поэтому количественное
поведение его градиентного спуска представляет отдельный интерес.

\subsection{Какой поток мы изучаем}

Пусть \(X\) имеет заданный исходный закон \(\mu\), а \(\Gamma\) ---
независимый от него стандартный гауссовский вектор.
Переменной задачи является случайный вектор \(Z\), измеримый
относительно \(\G=\sigma(X)\). Его рассматриваем в пространстве
\(\cH=L^2(\G;\R^d)\).
К \(Z\) добавляется гауссовский шум, после чего выбирается случайный
вектор \(Y\) с заданным целевым законом \(\nu\).
В поднятом функционале максимизируется среднее скалярное произведение
\(Z+\Gamma\) и \(Y\), а связь \(Z\) с исходным \(X\) учитывается
отдельным слагаемым:
\[
V(Z)=
\sup_{\substack{Y\in L^2(\sigma(X,\Gamma);\R^d)\\\Law(Y)=\nu}}
\left\{\E\ip{Z+\Gamma}{Y}-\E\ip{Z}{X}\right\}.
\]
Таким образом, нам важен не только закон \(Z\), но и его зависимость
от \(X\). Заменять второе слагаемое максимальной ковариацией между
их законами в общем случае нельзя: это изменило бы задачу.

Мы изучаем непрерывное уравнение
\[
\dot Z_t=-\nabla V(Z_t),\qquad Z_{t=0}=Z_0.
\]
Здесь \(t\) --- время движения к минимизатору функционала.
Оно не совпадает со временем мартингала между начальным и конечным
распределениями. Полная постановка, формула градиента и все условия
на исходные меры собраны в следующем разделе.

\subsection{От сходимости к скорости}

Пусть \(Z^\star\) --- предел траектории.
Мы хотим контролировать две ошибки: расстояние до этого предела
и оставшийся зазор функционала. Требуется найти \(\kappa>0\),
для которого одновременно при всех \(t\ge0\)
\[
\begin{aligned}
\norm{Z_t-Z^\star}_2
&\le\norm{Z_0-Z^\star}_2e^{-\kappa t},\\
V(Z_t)-\inf_{\cH}V
&\le\bigl(V(Z_0)-\inf_{\cH}V\bigr)e^{-\kappa t}.
\end{aligned}
\]
В обеих оценках множитель равен соответствующей начальной ошибке.
Оценка с произвольным дополнительным множителем или только для
достаточно больших времён отвечала бы на более слабый вопрос.

Мы получаем эти оценки для неприводимой пары мер, когда целевая
поддержка компактна и выпукла, плотность целевой меры почти всюду
на внутренности этой поддержки ограничена сверху и отделена от нуля,
а поддержка исходной меры лежит
строго внутри целевой. Начальное данное предполагается ограниченным.
Константа \(\kappa\) может зависеть от исходных данных и начального
данного. Её оптимальность и независимость от всех данных не утверждаются.

Условие на среднее тоже существенно. Из выпуклого порядка следует,
что у \(\mu\) и \(\nu\) одинаковые центры масс.
Поэтому прибавление постоянного вектора к \(Z\) не меняет значение
\(V\). Единственность минимизатора нужно понимать в классе с заданным
средним. При этом минимальное значение ищется по всему пространству \(\cH\),
включая неограниченные случайные векторы.

\subsection{Главная идея доказательства}

Для экспоненциальной скорости естественно искать сильную выпуклость.
В нашей задаче нужная оценка получается на множествах вида
\[
\{Z\in\cH:\ \norm{Z}_\infty\le R,\ \E Z=\bar z\}.
\]
Для каждого фиксированного радиуса \(R\) мы найдём положительную
константу кривизны. Известная качественная теорема даёт один радиус,
который ограничивает всю траекторию. После этого одну и ту же
константу кривизны можно использовать при всех временах.

Основная трудность состоит в получении самой количественной кривизны.
Мы временно вводим энтропию и для двух исходных точек берём
их оптимальные энтропийные потенциалы. Линейно интерполируя потенциалы
на концах отрезка, строим верхнее сравнение для значения функционала.
На концах значения совпадают; оптимальность промежуточного потенциала
нам не нужна.

Вторая производная этого сравнения выражается через дисперсию.
В ней присутствует неизвестная разность целевых потенциалов.
Сравнение двух независимых копий устраняет этот член.
Оставшаяся величина связана с гессианом логарифмического интеграла.
Геометрия целевой поддержки и нижняя граница её плотности дают
положительную интегральную оценку этого гессиана, а гауссовское
усреднение передаёт её в исходное пространство случайных векторов.

Важно, что полученная константа не зависит от энтропийного параметра.
Поэтому она сохраняется при удалении энтропии.
На ограниченной траектории сильная монотонность градиента затем даёт
обе требуемые оценки. Схема на рис.~\ref{ru:fig:bass-route}
показывает, где используется известная качественная сходимость,
а где возникает новый количественный шаг.

\begin{figure}[H]
\centering
\begin{tikzpicture}[
  font=\small,
  bass box/.style={draw=black!45,rounded corners=2pt,line width=.5pt,
    align=center,inner sep=7pt,minimum height=1.02cm},
  bass step/.style={bass box,text width=4.7cm,fill=geoPtwo!6},
  bass known/.style={bass box,text width=4.7cm,fill=black!3},
  bass arrow/.style={-{Stealth[length=2mm]},draw=black!65,line width=.65pt}
]
\node[bass box,text width=10.9cm,fill=black!2] (data) at (0,2.5)
  {Исходная пара мер и поднятый функционал \(V\)};
\node[bass step] (entropy) at (-3.05,.9)
  {Энтропийное сравнение\\и две независимые копии};
\node[bass known] (qual) at (3.05,.9)
  {Известный результат:\\поток, ограниченность,\\сильный предел};
\node[bass step] (curvature) at (-3.05,-.85)
  {Положительная кривизна \(m_R\)\\при \(\norm{Z}_\infty\le R\), \(\E Z=\bar z\)};
\node[bass known] (radius) at (3.05,-.85)
  {Один радиус для всей орбиты:\\\(\norm{Z_t}_\infty\le R\), \(t\ge0\)};
\node[bass box,text width=10.9cm,fill=geoPtwo!9] (rate) at (0,-2.8)
  {Сильная монотонность вдоль орбиты\\
   \(L^2\)-ошибка и зазор функционала убывают экспоненциально};
\draw[bass arrow] (data.south) -- ++(0,-.24) -| (entropy.north);
\draw[bass arrow] (data.south) -- ++(0,-.24) -| (qual.north);
\draw[bass arrow] (entropy) -- (curvature);
\draw[bass arrow] (qual) -- (radius);
\draw[bass arrow] (curvature.south) -- ++(0,-.32) -| ([xshift=-2.4cm]rate.north);
\draw[bass arrow] (radius.south) -- ++(0,-.32) -| ([xshift=2.4cm]rate.north);
\end{tikzpicture}
\caption{Схема доказательства скорости. Слева показан количественный аргумент
настоящей работы, справа --- используемый известный качественный результат.
Стрелки обозначают зависимости: кривизна применяется с радиусом,
который ограничивает всю орбиту.}
\label{ru:fig:bass-route}
\end{figure}

\subsection{Устройство работы}

В разделе~2 мы точно формулируем задачу и основной результат,
а также выделяем используемые свойства известного потока.
В разделе~3 вводим энтропийное сравнение и обосновываем переход
к исходному функционалу. В разделе~4 получаем количественную кривизну,
включая оценку двумя копиями. В разделе~5 выводим экспоненциальную
сходимость. Раздел~6 посвящён актуальности результата и его связи
с другими работами. В разделе~7 описаны формальная верификация
и повторная проверка исходников, в разделе~8 --- использование ИИ.
Раздел~9 завершает работу.

\section{Постановка задачи и основной результат}\label{ru:sec:setting}

\subsection{Какие данные мы фиксируем}

Пусть $d\ge2$, а $\mu,\nu\in\mathcal P_2(\R^d)$ --- вероятностные меры
с конечными вторыми моментами. Будем считать, что выполнены следующие условия.
\begin{enumerate}[label=(A\arabic*),leftmargin=2.7em]
\item $\mu\le_{cx}\nu$: для каждой выпуклой функции, для которой определены
соответствующие интегралы, интеграл по $\mu$ не превосходит интеграла по $\nu$.
\item Пара $(\mu,\nu)$ неприводима: для любых измеримых $A,B\subset\R^d$
с $\mu(A)>0$ и $\nu(B)>0$ существует мартингальная связка $\pi$ этих мер,
для которой $\pi(A\times B)>0$.
Мартингальность означает $\int y\,\pi_x(dy)=x$ для $\mu$-почти всех $x$.
Связка в этом условии может зависеть от обоих множеств.
\item Множество $C=\supp\nu$ компактно и выпукло;
$\supp\mu\subset I:=\operatorname{int}C$.
\item Мера $\nu$ абсолютно непрерывна относительно $d$-мерной меры Лебега.
Для её плотности на $I$ существуют константы $0<\underline f\le\overline f<\infty$,
такие что $\underline f\le f_\nu\le\overline f$ почти всюду.
\end{enumerate}

Последнее условие, в частности, означает, что $C$ полноразмерно.
Граница выпуклого тела имеет нулевую меру Лебега, поэтому $\nu(I)=1$.
Из компактности $\supp\mu$ и его внутреннего положения следует положительный
зазор до границы $C$. Гладкости этой границы или плотности мы не предполагаем.
Условия соответствуют многомерному продолжению теоремы~1.9
из~\citep{ru:bps2025}. Неприводимость здесь записана через связки;
в записи через мартингал достаточно взять его совместный закон в моменты $0$ и $1$.
Обратно, связка реализуется двухэтапным мартингалом с теми же концами.

На произвольном вероятностном пространстве фиксируем независимые векторы
$X\sim\mu$ и $\Gamma\sim\gamma=N(0,I_d)$. Положим
\[
 \G=\sigma(X),\qquad \mathcal F=\sigma(X,\Gamma),\qquad
 \cH=L^2(\G;\R^d).
\]
Все равенства между случайными величинами понимаются с точностью до почти
верного совпадения. Скалярное произведение в $\cH$ равно
$\ip{Z}{W}_{\cH}=\E\ip{Z}{W}$; его норму обозначаем $\norm{Z}_2$.

\subsection{Почему нужен именно поднятый функционал}

Для $Z\in\cH$ определим
\begin{equation}\label{ru:eq:literal-functional}
 V(Z)=\sup_{\substack{Y\in L^2(\mathcal F;\R^d)\\\Law(Y)=\nu}}
 \E\bigl[\ip{Z+\Gamma}{Y}-\ip{Z}{X}\bigr].
\end{equation}
Первое слагаемое ищет наиболее выгодное сопряжение гауссовски сглаженного
наблюдения с целевым законом. Во втором слагаемом векторы $Z$ и $X$ остаются
на исходном вероятностном пространстве. Именно эту связь надо сохранить:
она входит в функционал, а не выбирается заново при каждом вычислении.

Если $\MCov(\rho,\nu)=\sup_{\pi\in\Pi(\rho,\nu)}\int\ip{w}{y}\,d\pi$,
то достижение оптимального транспорта даёт
\begin{equation}\label{ru:eq:mcov-functional}
 V(Z)=\MCov(\Law(Z)*\gamma,\nu)-\E\ip{Z}{X}.
\end{equation}
Действительно, $Z$ измерим относительно $\sigma(X)$, поэтому независимость
обеспечивает $\Law(Z+\Gamma)=\Law(Z)*\gamma$. Этот закон имеет всюду положительную
плотность, полученную свёрткой с гауссовской мерой. Оптимальная карта Бренье $T_Z$ из него в $\nu$ задаёт
допустимый в~\eqref{ru:eq:literal-functional} вектор $Y=T_Z(Z+\Gamma)$ и достигает
максимальной ковариации. Так совпадают оба значения супремума.

Из~\eqref{ru:eq:mcov-functional} не следует, что $V(Z)$ зависит только от закона $Z$:
слагаемое $\E\ip{Z}{X}$ по-прежнему учитывает его связь с $X$.
В частности, мы не подменяем его величиной $\MCov(\Law(Z),\mu)$.

Градиент известен явно~\citep[лемма~2.1]{ru:bps2025}:
\begin{equation}\label{ru:eq:actual-gradient}
 g(Z):=\nabla V(Z)=\E[T_Z(Z+\Gamma)\mid\G]-X.
\end{equation}
Это непрерывный градиент Фреше на всём $\cH$, а $V$ выпуклый.
Чтобы увидеть роль этой формулы, зафиксируем оптимальный вектор $Y_Z$.
Он остаётся допустимым при замене $Z$ на $W$, откуда
\begin{equation}\label{ru:eq:global-support}
 V(W)\ge V(Z)+\ip{g(Z)}{W-Z}_{\cH}\qquad(Z,W\in\cH).
\end{equation}
Два таких сравнения дают точную оценку остатка:
\[
0\le V(W)-V(Z)-\ip{g(Z)}{W-Z}_{\cH}
\le\norm{g(W)-g(Z)}_2\norm{W-Z}_2.
\] Непрерывность выбранного оптимального
вектора, а значит и $g$, позволяет получить производную без дифференцирования
оптимизатора по параметру. При компактной цели эта непрерывность следует из
устойчивости оптимальных планов и единственности карты Бренье.

Из выпуклого порядка следуют равные барицентры $\E X=\int y\,\nu(dy)$:
достаточно применить его к обоим знакам линейной функции. Поэтому
\[
 V(Z+c)=V(Z),\qquad \E g(Z)=0\qquad(c\in\R^d).
\]
Следовательно, сильной выпуклости на всём $\cH$ быть не может: вдоль
постоянных переносов функционал не меняется. Нужная количественная кривизна
возникает после фиксации среднего. Для радиуса $R$ и вектора $\bar z$ обозначим
\[
 S_{R,\bar z}=\{Z\in\cH:\norm{Z}_\infty\le R,\ \E Z=\bar z\}.
\]
Это выпуклое замкнутое подмножество $\cH$. Замкнутость условия $\norm{Z}_\infty\le R$
следует из почти всюду сходящейся подпоследовательности любой
сходящейся в $L^2$ последовательности.

\subsection{Основная теорема}

Рассматриваем градиентное время $t\in[0,\infty)$ и уравнение
\begin{equation}\label{ru:eq:flow}
 \dot Z_t=-g(Z_t),\qquad Z_{t=0}=Z_0.
\end{equation}
Это время движения по функционалу. Оно не совпадает с временем мартингала
Басса, которое обычно меняется от $0$ до $1$.

\begin{theoremru}[Многомерная экспоненциальная сходимость]\label{ru:thm:main}
При условиях (A1)--(A4) для каждого $Z_0\in\cH\cap L^\infty$
существует единственное локально абсолютно непрерывное решение
\eqref{ru:eq:flow}, удовлетворяющее уравнению почти всюду по времени.
Оно является $C^1$-кривой в $\cH$, сохраняет среднее и равномерно ограничено
в $L^\infty$ на всей полуоси. Существует глобальный минимизатор $Z^\star$ функционала
$V$ на всём $\cH$, единственный среди глобальных минимизаторов со средним
$\E Z_0$, и константа $\kappa>0$, для которых при всех $t\ge0$
\begin{align}
 \norm{Z_t-Z^\star}_2
 &\le\norm{Z_0-Z^\star}_2 e^{-\kappa t},\label{ru:eq:main-norm}\\
 V(Z_t)-\inf_{W\in\cH}V(W)
 &\le\left[V(Z_0)-\inf_{W\in\cH}V(W)\right]e^{-\kappa t}.
 \label{ru:eq:main-energy}
\end{align}
Константа может зависеть от исходных данных и начального состояния.
\end{theoremru}

В обеих оценках множителем служит именно начальная ошибка. Утверждение
не ограничивается достаточно большими временами. Начальное среднее произвольно;
монотонности сопряжения $Z_0$ с $X$, конечности поддержки $\mu$ или малости
$\norm{Z_0}_\infty$ не требуется.

\paragraph{На какой известный результат мы опираемся.}
Теорема~1.8 из~\citep{ru:bps2025}, применённая к тому же функционалу,
уже даёт существование и единственность потока, его сильный $L^2$-предел,
глобальную минимальность и единственность минимизатора при фиксированном среднем.
Для ограниченного старта она также даёт один существенный радиус для всей
траектории. Непрерывность~\eqref{ru:eq:actual-gradient} улучшает локально
абсолютно непрерывную кривую до $C^1$; равенство $\E g=0$ сохраняет среднее.
Единственность среди локально абсолютно непрерывных решений можно также
проверить непосредственно: монотонность градиента выпуклого функционала
делает квадрат расстояния между двумя такими решениями невозрастающим
почти всюду; при одинаковом старте он равен нулю.
Эту качественную часть мы используем как опубликованную теорему.
Новым шагом в дальнейшем доказательстве будет положительная кривизна на
$S_{R,\bar z}$ и вытекающие из неё скорости.

В~\citep{ru:bps2025} меры предварительно центрированы. Это не ограничивает
применимость результата: при общем барицентре $m$ заменим $X,Y$ на
$X-m,Y-m$ и перенесём обе меры на $-m$. Условия сохраняются, а значение
поднятого функционала не меняется, поскольку
\[
 \E\bigl[\ip{Z+\Gamma}{Y-m}-\ip{Z}{X-m}\bigr]
 =\E\bigl[\ip{Z+\Gamma}{Y}-\ip{Z}{X}\bigr]-\ip{\E\Gamma}{m}
 =\E\bigl[\ip{Z+\Gamma}{Y}-\ip{Z}{X}\bigr].
\]
При этом $\sigma(X-m)=\sigma(X)$, так что пространство $\cH$ также
остаётся тем же.

Такой порядок не предполагает заранее ни скорость, ни сильную выпуклость.
Из известной теоремы сначала получаем радиус $R$, включающий поток и его
предел, затем доказываем применимую на этом радиусе количественную оценку.
Раздел~\ref{ru:sec:rates} завершит доказательство теоремы~\ref{ru:thm:main}.

\section{Энтропийное сравнение и его концевые точки}\label{ru:sec:ent:comparison}

Чтобы получить скорость сходимости, нам понадобится положительная нижняя
оценка кривизны функционала. Но карта оптимального переноса определяется
вариационной задачей, и непосредственно дифференцировать её по исходному
случайному вектору неудобно. Начнём поэтому с регуляризованной задачи, в которой
нужная производная записывается через обычную дисперсию.

Обозначим через \(D(\pi\Vert\rho\otimes\nu)\) относительную энтропию
сопряжения \(\pi\) относительно произведения его маргиналей. При
\(\varepsilon>0\) положим
\[
 E_\varepsilon(\rho,\nu)=
 \sup_{\substack{\pi\in\Pi(\rho,\nu)\\
                  D(\pi\Vert\rho\otimes\nu)<\infty}}
 \left\{\int\ip{w}{y}\,\pi(dw,dy)
       -\varepsilon D(\pi\Vert\rho\otimes\nu)\right\}.
\]
Здесь \(\Pi(\rho,\nu)\) --- множество всех сопряжений мер \(\rho\) и \(\nu\).
Поскольку целевой вектор ограничен, для конечности ковариационного слагаемого
достаточно первого момента \(\rho\).

\subsection{Строка, столбец и два маргинальных уравнения}

Пусть \(C=\operatorname{supp}\nu\subset\overline B_M(0)\).
Целевой потенциал позволяет записать верхнюю оценку значения задачи.
Сначала зафиксируем его произвольно, а затем выберем так, чтобы оценка
достигалась в нужных концевых состояниях.
Для ограниченной измеримой функции \(a\colon C\to\R\) введём
\[
 v_a(w)=\varepsilon\log\int_C
 e^{(\ip{w}{y}-a(y))/\varepsilon}\,\nu(dy),
 \qquad
 p_{a,w}(dy)=
 e^{(\ip{w}{y}-a(y)-v_a(w))/\varepsilon}\,\nu(dy).
\]
Функцию \(a\) будем называть столбцом, а \(v_a\) --- соответствующей строкой.
По определению \(p_{a,w}\) является вероятностной мерой.
Добавление константы \(c\) к столбцу вычитает эту константу из строки и
сохраняет плотность сопряжения. Этой свободой будем пользоваться, чтобы
нормировать строку условием \(v_a(0)=0\).

Для любого такого столбца и любой исходной меры $\rho$ с первым моментом
мера $Q_a(dw,dy)=\rho(dw)p_{a,w}(dy)$ является вероятностной.
Её первая маргиналь равна $\rho$; вторая пока может отличаться от $\nu$.
Для связки $\pi\in\Pi(\rho,\nu)$ с конечной относительной энтропией
неотрицательность $D(\pi\Vert Q_a)$ даёт слабую двойственность:
\begin{equation}\label{ru:eq:ent:weak-duality}
 \int\ip{w}{y}\,d\pi-\varepsilon D(\pi\Vert\rho\otimes\nu)
 \le\int v_a\,d\rho+\int a\,d\nu.
\end{equation}
Действительно, логарифм плотности $Q_a$ относительно $\rho\otimes\nu$
равен $(\ip{w}{y}-a(y)-v_a(w))/\varepsilon$;
подстановка этой формулы в $D(\pi\Vert Q_a)\ge0$ даёт указанное неравенство.
Все слагаемые интегрируемы благодаря ограниченности цели и столбца,
первому моменту $\rho$ и липшицевости строки.
Именно поэтому верхнее сравнение допустимо и для промежуточного
неоптимального столбца.

\begin{lemmaru}\label{ru:lem:ent:calculus}
Функция \(v_a\) дважды непрерывно дифференцируема, выпукла и \(M\)-липшицева. Для любых \(w,q,q'\in\R^d\)
\begin{align}
 Dv_a(w)[q]&=\int\ip{q}{y}\,p_{a,w}(dy),\label{ru:eq:ent:first-derivative}\\
 D^2v_a(w)[q,q']&=\varepsilon^{-1}
 \operatorname{Cov}_{p_{a,w}}\bigl(\ip{q}{y},\ip{q'}{y}\bigr).
 \label{ru:eq:ent:hessian}
\end{align}
\end{lemmaru}
\begin{proof}
На каждом ограниченном множестве значений \(w\) подынтегральная экспонента
и все её производные ограничены общей константой, умноженной на степень
\(\norm{y}\). Поскольку \(\norm{y}\le M\) и \(\nu\) --- вероятность, эти
мажоранты интегрируемы. Можно последовательно дифференцировать интеграл.
Производная логарифма даёт первую формулу; ещё одно дифференцирование
даёт разность второго момента и произведения первых моментов, то есть
ковариацию. При \(q'=q\) она неотрицательна, поэтому строка выпукла.
Наконец, её градиент есть среднее векторов из \(\overline B_M(0)\),
так что его норма не превосходит \(M\).
\end{proof}

Нам нужны не только такие формальные строки, но и пары, действительно
решающие регуляризованную задачу. Следующее утверждение обеспечивает
равенство в двух концах отрезка, который позднее будем интерполировать.

\begin{propositionru}\label{ru:prop:ent:endpoints}
Пусть \(Z\) --- ограниченный случайный вектор, независимый от стандартного
гауссовского вектора \(\Gamma\), \(W=Z+\Gamma\) и \(\rho=\Law(W)\).
Для каждого \(\varepsilon>0\) существует столбец \(a\), ограниченный и
непрерывный на \(C\), выпуклый на \(C\), для которого
\begin{align}
 v_a(0)&=0,\label{ru:eq:ent:normalization}\\
 \E e^{(\ip{W}{y}-v_a(W))/\varepsilon}
   &=e^{a(y)/\varepsilon}\qquad(y\in C).\label{ru:eq:ent:column}
\end{align}
Мера
\[
 \pi_a(dw,dy)=
 e^{(\ip{w}{y}-a(y)-v_a(w))/\varepsilon}\,\rho(dw)\nu(dy)
\]
имеет маргинали \(\rho,\nu\), конечную относительную энтропию и достигает
значения \(E_\varepsilon(\rho,\nu)\). В частности,
\begin{equation}\label{ru:eq:ent:primal-equality}
 E_\varepsilon(\rho,\nu)=\E v_a(W)+\int a\,d\nu.
\end{equation}
\end{propositionru}

\begin{proof}
Сначала рассмотрим конечный источник
\(\rho=\sum_{i=1}^n p_i\delta_{w_i}\), где \(p_i>0\). Минимизируем функцию
\[
 G(b)=\sum_i p_i b_i+
 \varepsilon\int_C\log\left(
 \sum_i p_i e^{(\ip{w_i}{y}-b_i)/\varepsilon}\right)\nu(dy)
 \quad(b\in\R^n).
\]
Она не меняется при добавлении общей константы ко всем \(b_i\).
Выберем нормировку \(\min_i b_i=0\). Если
\(\lvert\ip{w_i}{y}\rvert\le K\) на \(C\) и
\(p_{\min}=\min_i p_i\), то
\[
 G(b)\ge p_{\min}\max_i b_i-K+\varepsilon\log p_{\min}.
\]
Действительно, первое слагаемое не меньше
\(p_{\min}\max_i b_i\), а в сумме под логарифмом можно оставить индекс,
для которого \(b_i=0\). Поэтому замкнутые подуровни в выбранной нормировке
компактны и минимум достигается. Инвариантность при общих сдвигах означает,
что найденный вектор является минимизатором и без ограничения нормировки.

Положим
\[
 a(y)=\varepsilon\log\sum_i p_i
 e^{(\ip{w_i}{y}-b_i)/\varepsilon}.
\]
Условия первого порядка для каждого \(b_i\) дают
\[
 \int_C e^{(\ip{w_i}{y}-a(y))/\varepsilon}\,\nu(dy)
 =e^{b_i/\varepsilon}.
\]
Следовательно, \(v_a(w_i)=b_i\). Определение \(a\) даёт второе
маргинальное уравнение. Общим сдвигом пары нормируем \(v_a(0)=0\).
Функция \(a\) выпукла и непрерывна, поскольку является логарифмом
интеграла экспонент аффинных функций.

Теперь вернёмся к \(W=Z+\Gamma\). Выберем конечные простые приближения
\(W_n\to W\) почти наверное и в \(L^2\), удовлетворяющие
\(\norm{W_n}\le2\norm{W}\). Для каждого приближения построим пару
\((a_n,v_n)\) с \(v_n(0)=0\). Из леммы~\ref{ru:lem:ent:calculus}
следует \(\lvert v_n(w)\rvert\le M\norm{w}\).

Покажем, почему столбцы имеют общую оценку и общую непрерывность при
фиксированном \(\varepsilon\). Их второе маргинальное уравнение имеет вид
\[
 a_n(y)=\varepsilon\log
 \E e^{(\ip{W_n}{y}-v_n(W_n))/\varepsilon}.
\]
На \(C\) модуль выражения в экспоненте не больше
\(2M\norm{W_n}/\varepsilon\le4M\norm{W}/\varepsilon\).
Вектор \(W\) имеет все экспоненциальные моменты нормы, поскольку \(Z\)
ограничен, а \(\Gamma\) гауссовский. Верхняя оценка для интеграла поэтому
едина по \(n,y\). Нижнюю оценку даёт Йенсен:
\[
 a_n(y)\ge-2M\E\norm{W_n}\ge-4M\E\norm{W}.
\]
При изменении \(y\) производная интеграла сверху контролируется
\(\E[\norm{W_n}e^{2M\norm{W_n}/\varepsilon}]/\varepsilon\).
В сочетании с положительной нижней оценкой самого интеграла это даёт
общую липшицеву оценку для \(a_n\) на \(C\). Все использованные моменты
конечны и ограничены независимо от \(n\).

По теореме Асколи выберем подпоследовательность, для которой
\(a_n\to a\) равномерно на \(C\). Из определения строки
\[
 \sup_{w\in\R^d}\lvert v_{a_n}(w)-v_a(w)\rvert
 \le\sup_{y\in C}\lvert a_n(y)-a(y)\rvert.
\]
Поэтому \(v_a(0)=0\), а
\(v_n(W_n)\to v_a(W)\) почти наверное. В столбцовом уравнении
можно перейти к пределу по той же гауссовской экспоненциальной мажоранте.
Получаем~\eqref{ru:eq:ent:column} для каждого \(y\in C\).
Предельный столбец выпуклый и непрерывный на \(C\).

Первое маргинальное уравнение выполнено по определению строки,
второе --- по~\eqref{ru:eq:ent:column}. Значит \(\pi_a\) имеет нужные
маргинали. Логарифм её плотности интегрируем относительно \(\pi_a\):
столбец ограничен, \(\lvert v_a(w)\rvert\le M\norm{w}\), а
\(\lvert\ip{w}{y}\rvert\le M\norm{w}\). Поэтому энтропия конечна.

Осталось проверить оптимальность именно этой пары. Для любого
сопряжения \(\pi\) с конечной относительной энтропией неотрицательность
\(D(\pi\Vert\pi_a)\) даёт
\[
 \int\ip{w}{y}\,d\pi-\varepsilon D(\pi\Vert\rho\otimes\nu)
 \le \int v_a\,d\rho+\int a\,d\nu .
\]
Все аффинные и строковые слагаемые интегрируемы; поэтому вычитание
логарифма плотности корректно. При \(\pi=\pi_a\) получается равенство.
Это доказывает и достижение, и~\eqref{ru:eq:ent:primal-equality}.
\end{proof}

\begin{remarkru}\label{ru:rem:ent:same-endpoints}
В дальнейшем мы будем менять столбец между двумя такими концевыми
значениями. Промежуточный столбец может не быть оптимальным. Нам понадобится
только слабая двойственность внутри отрезка и равенство на его концах.
Поэтому дифференцировать семейство оптимизаторов не потребуется.
\end{remarkru}

\subsection{Сравнивающая функция}\label{ru:sec:ent:scalar}

Зафиксируем ограниченные $Z_0,Z_1\in\cH$ и $\varepsilon>0$.
Для источников $\Law(Z_i+\Gamma)$ выберем столбцы $a_i$ из
утверждения~\ref{ru:prop:ent:endpoints} и положим
$a_s=(1-s)a_0+sa_1$, $0\le s\le1$.

Положим \(Z_s=(1-s)Z_0+sZ_1\), \(h=Z_1-Z_0\), \(b=a_1-a_0\) и рассмотрим
\[
 J(s)=\E v_{a_s}(Z_s+\Gamma)+\int_C a_s\,d\nu .
\]
По слабой двойственности
\(E_\varepsilon(\Law(Z_s+\Gamma),\nu)\le J(s)\), а при \(s=0,1\)
выполнено равенство. Поскольку \(h\), \(b\) и целевой вектор ограничены,
можно дважды дифференцировать по \(s\) под знаком ожидания. Например,
первая производная строки контролируется
\(M\norm{h}_\infty+\norm{b}_\infty\), а вторая --- квадратом
этой константы, делённым на \(\varepsilon\). Получаем
\begin{equation}\label{ru:eq:curv:comparison-second}
 J''(s)=\frac1\varepsilon\E
 \operatorname{Var}_{p_{a_s,Z_s+\Gamma}}
       \bigl(\ip{h}{y}-b(y)\bigr).
\end{equation}

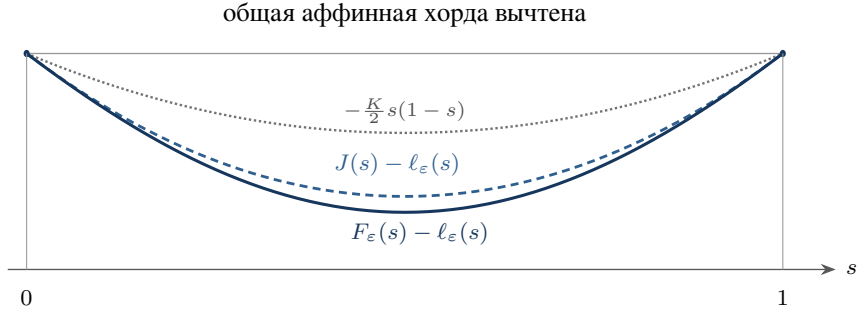
\begin{figure}[H]
\centering
\begin{tikzpicture}[x=10cm,y=4.2cm,font=\small]
  \draw[-{Stealth[length=2mm]},draw=black!65,line width=.5pt]
    (-.025,-.68) -- (1.07,-.68) node[right] {\(s\)};
  \draw[draw=black!40,line width=.5pt] (0,0) -- (1,0)
    node[midway,above=6pt] {общая аффинная хорда вычтена};
  \draw[draw=black!55,densely dotted,line width=.9pt]
    plot[domain=0:1,samples=65] (\x,{-\x*(1-\x)});
  \draw[draw=geoPthreehalf,densely dashed,line width=1.05pt]
    plot[domain=0:1,samples=65] (\x,{-1.8*\x*(1-\x)});
  \draw[draw=geoPtwo,line width=1.15pt]
    plot[domain=0:1,samples=65]
      (\x,{-1.8*\x*(1-\x)-.8*\x*\x*(1-\x)*(1-\x)});
  \fill[geoPtwo] (0,0) circle[x radius=.004,y radius=.0114];
  \fill[geoPtwo] (1,0) circle[x radius=.004,y radius=.0114];
  \node[below=4pt] at (0,-.68) {\(0\)};
  \node[below=4pt] at (1,-.68) {\(1\)};
  \node[align=center,text=black!70,font=\footnotesize] at (.5,-.18)
    {\(-\frac K2s(1-s)\)};
  \node[align=center,text=geoPthreehalf] at (.49,-.35)
    {\(J(s)-\ell_\varepsilon(s)\)};
  \node[align=center,text=geoPtwo] at (.52,-.565)
    {\(F_\varepsilon(s)-\ell_\varepsilon(s)\)};
  \draw[draw=black!40,line width=.4pt] (0,0) -- (0,-.68);
  \draw[draw=black!40,line width=.4pt] (1,0) -- (1,-.68);
\end{tikzpicture}
\caption{Схематическое сравнение при фиксированном \(\varepsilon>0\).
Здесь \(F_\varepsilon(s)=E_\varepsilon(\Law(Z_s+\Gamma),\nu)\),
а \(J(s)\) --- верхнее сравнение, построенное по интерполированным
целевым потенциалам. Вычтена общая хорда
\(\ell_\varepsilon(s)=(1-s)F_\varepsilon(0)+sF_\varepsilon(1)\).
Если \(J''\ge K\), где \(K>0\), то
\(F_\varepsilon(s)-\ell_\varepsilon(s)
\le J(s)-\ell_\varepsilon(s)\le-(K/2)s(1-s)\).
На концах все три выражения равны нулю.
Кривые поясняют направления неравенств и не являются вычисленными
значениями энтропийного функционала.}
\label{ru:fig:bass-comparison}
\end{figure}

\subsection{Удаление энтропии}

Положительная нижняя оценка кривизны должна сохраниться при
\(\varepsilon\downarrow0\). Сначала проверим сам переход к
нерегуляризованной ковариационной задаче.

\begin{propositionru}\label{ru:prop:ent:removal}
Для любой вероятности \(\rho\) с конечным первым моментом и компактно
поддержанной \(\nu\)
\[
 E_\varepsilon(\rho,\nu)\longrightarrow\MCov(\rho,\nu)
 \qquad(\varepsilon\downarrow0).
\]
\end{propositionru}
\begin{proof}
Верхняя оценка \(E_\varepsilon\le\MCov\) следует из неотрицательности
энтропии. Для обратной оценки зафиксируем \(a>0\) и разобьём \(C\)
на конечное число измеримых ячеек \(C_j\) диаметра не больше \(a\).
Ячейки нулевой массы отбросим. Для произвольного сопряжения
\(\pi\in\Pi(\rho,\nu)\) определим
\[
 \pi^a=\sum_{\nu(C_j)>0}
 \pi(\,\cdot\,\times C_j)\otimes\frac{\nu|_{C_j}}{\nu(C_j)}.
\]
Обе маргинали сохраняются. Если
\(r_j=d\pi(\,\cdot\,\times C_j)/d\rho\), то \(0\le r_j\le1\).
Поэтому плотность \(\pi^a\) относительно \(\rho\otimes\nu\)
на \(C_j\) равна \(r_j(w)/\nu(C_j)\) и не превосходит
\(1/\min_{\nu(C_j)>0}\nu(C_j)\). Её энтропия ограничена конечным
числом \(H_a\), не зависящим от исходного сопряжения.

При пересэмплировании обе целевые точки остаются в одной ячейке,
так что изменение \(\ip{w}{y}\) по модулю не больше \(a\norm{w}\).
Следовательно,
\[
 \int\ip{w}{y}\,d\pi^a
 \ge\int\ip{w}{y}\,d\pi-a\int\norm{w}\,\rho(dw).
\]
Беря супремум по \(\pi\), получаем
\[
 \MCov(\rho,\nu)-a\int\norm{w}\,d\rho-\varepsilon H_a
 \le E_\varepsilon(\rho,\nu)\le\MCov(\rho,\nu).
\]
Сначала устремим \(\varepsilon\) к нулю при фиксированном \(a\),
а затем \(a\) к нулю. Равномерной оценки для \(H_a\) при \(a\downarrow0\)
не требуется.
\end{proof}

Таким образом, мы можем работать с гладкими строками при положительном
\(\varepsilon\), если найдём оценку, которая не вырождается при его
стремлении к нулю. Именно такую оценку получим в следующем разделе.

\section{Положительная кривизна на ограниченном сечении}\label{ru:sec:curv:main}

Зафиксируем радиус \(0\le R<\infty\). Мы хотим сравнить два исходных вектора
\(Z_0,Z_1\in\cH\), для которых \(\norm{Z_i}_\infty\le R\).
Нижняя оценка должна зависеть от этого общего радиуса, но не от выбора
самих векторов. Сначала получим геометрическую оценку для строк, а затем
перенесём её на случайные исходные данные.

\subsection{Внутренняя оценка столбца и рост строки}

Выберем \(y_c\in I\) и \(r>0\) так, чтобы
\(\overline B_{2r}(y_c)\subset I\). Плотность \(\nu\) на этом шаре
не меньше некоторого \(m>0\). Пусть
\(G_d=\E\norm{\Gamma}\) и \(C\subset\overline B_M(0)\).

\begin{lemmaru}\label{ru:lem:curv:cap-growth}
Для концевых столбцов \(a_0,a_1\) из
утверждения~\ref{ru:prop:ent:endpoints}, отвечающих \(Z_0+\Gamma,Z_1+\Gamma\),
и их интерполяции \(a_s=(1-s)a_0+sa_1\), \(0\le s\le1\), существуют
константы \(c>0,K<\infty\), независимые от \(s\), от концов и от
\(0<\varepsilon\le1\), такие что
\begin{equation}\label{ru:eq:curv:row-growth}
 v_{a_s}(w)\ge\ip{w}{y_c}+c\norm{w}-K,\qquad v_{a_s}(0)\le0 .
\end{equation}
Все эти строки выпуклы и \(M\)-липшицевы.
\end{lemmaru}
\begin{proof}
Положим \(K_0=2M(R+G_d)\). Столбцовое уравнение и Йенсен дают
\[
 a_i(y)\ge\E[\ip{Z_i+\Gamma}{y}-v_{a_i}(Z_i+\Gamma)]\ge-K_0.
\]
Из концевого равенства, оценки
\(E_\varepsilon(\Law(Z_i+\Gamma),\nu)\le M\E\norm{Z_i+\Gamma}\)
и \(\lvert v_{a_i}(w)\rvert\le M\norm{w}\) следует
\(\int a_i\,d\nu\le K_0\). Эти две оценки сохраняются при
интерполяции столбцов.

Для \(y\in\overline B_r(y_c)\) шар \(B_r(y)\) лежит в
\(B_{2r}(y_c)\). Функция \(a_s+K_0\) выпукла и неотрицательна.
По Йенсену для равномерной меры на шаре и нижней оценке плотности
\[
 a_s(y)+K_0\le\frac1{|B_r|}\int_{B_r(y)}(a_s+K_0)\,dy
 \le\frac1{m|B_r|}\int_C(a_s+K_0)\,d\nu
 \le\frac{2K_0}{m|B_r|}.
\]
Обозначим правую верхнюю оценку для \(a_s(y)\) через \(A\).
При \(w\ne0\) шар радиуса \(r/4\) с центром
\(y_c+(r/2)w/\norm{w}\) лежит в \(B_r(y_c)\), а на нём
\(\ip{w}{y}\ge\ip{w}{y_c}+(r/4)\norm{w}\).
Его масса не меньше \(q=m|B_{r/4}|>0\), причём \(q\le1\).
Оставляя этот шар в строковом интеграле, получаем
\[
 v_{a_s}(w)\ge\ip{w}{y_c}+(r/4)\norm{w}-A+\varepsilon\log q
 \ge\ip{w}{y_c}+(r/4)\norm{w}-A+\log q .
\]
При \(w=0\) та же оценка следует непосредственно из интеграла.
Можно взять \(c=r/4\), \(K=A-\log q\).

Наконец, концевые нормировки означают
\(\int e^{-a_i/\varepsilon}\,d\nu=1\). Выпуклость экспоненты даёт
\[
 \int e^{-a_s/\varepsilon}\,d\nu
 \le(1-s)\int e^{-a_0/\varepsilon}\,d\nu
       +s\int e^{-a_1/\varepsilon}\,d\nu=1.
\]
Поэтому \(v_{a_s}(0)\le0\). Выпуклость и липшицевость следуют из
леммы~\ref{ru:lem:ent:calculus}.
\end{proof}

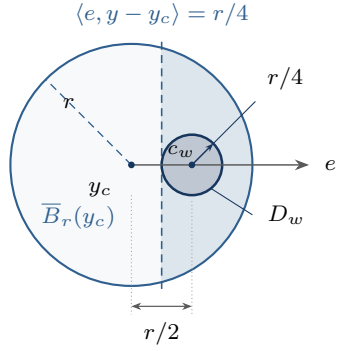
\begin{figure}[H]
\centering
\begin{tikzpicture}[font=\small,
  dim/.style={{Stealth[length=1.3mm]}-{Stealth[length=1.3mm]},draw=black!65,line width=.45pt}]
  \fill[geoPthreehalf!4] (0,0) circle[radius=1.6];
  \begin{scope}
    \clip (0,0) circle[radius=1.6];
    \fill[geoPthreehalf!16] (.4,-1.65) rectangle (1.65,1.65);
  \end{scope}
  \draw[draw=geoPthreehalf,line width=.8pt] (0,0) circle[radius=1.6];
  \draw[draw=geoPthreehalf,densely dashed,line width=.6pt] (.4,-1.64) -- (.4,1.64);
  \node[anchor=south,text=geoPthreehalf] at (.4,1.71)
    {\(\ip{e}{y-y_c}=r/4\)};
  \filldraw[draw=geoPtwo,fill=geoPtwo!28,line width=.95pt]
    (.8,0) circle[radius=.4];
  \draw[-{Stealth[length=1.9mm]},draw=black!65,line width=.55pt]
    (0,0) -- (2.35,0) node[right=2pt] {\(e\)};
  \draw[draw=geoPthreehalf,densely dashed,line width=.5pt]
    (0,0) -- (-1.13137,1.13137) node[midway,above left=2pt] {\(r\)};
  \fill[geoPtwo] (0,0) circle[radius=1.25pt];
  \node[anchor=north east] at (-.12,-.10) {\(y_c\)};
  \fill[geoPtwo] (.8,0) circle[radius=1.25pt];
  \node[font=\footnotesize] at (.66,.18) {\(c_w\)};
  \node[text=geoPthreehalf] at (-.70,-.67) {\(\overline B_r(y_c)\)};
  \draw[-{Stealth[length=1.1mm]},draw=geoPtwo,line width=.55pt]
    (.8,0) -- (1.08284,.28284);
  \draw[draw=geoPtwo,line width=.4pt] (1.08284,.28284) -- (1.64,.84)
    node[above right=1pt] {\(r/4\)};
  \draw[draw=geoPtwo,line width=.4pt] (1.05,-.31) -- (1.60,-.66)
    node[right=2pt] {\(D_w\)};
  \draw[draw=black!35,densely dotted,line width=.4pt] (0,-.40) -- (0,-1.98);
  \draw[draw=black!35,densely dotted,line width=.4pt] (.8,-.48) -- (.8,-1.98);
  \draw[dim] (0,-1.87) -- (.8,-1.87);
  \node[below=2pt] at (.4,-1.87) {\(r/2\)};
  \node[anchor=west,align=left,text width=6.7cm] at (3.05,1.13)
    {\(e=w/\norm{w},\quad w\ne0\),\\[4pt]
     \(c_w=y_c+(r/2)e,\quad D_w=B_{r/4}(c_w)\).};
  \node[anchor=west,align=left,text width=6.7cm] at (3.05,-.12)
    {Для каждой точки \(y\in D_w\):\\[4pt]
     \(\ip{w}{y-y_c}\ge(r/4)\norm{w}\).};
  \node[anchor=west,align=left,text width=6.7cm] at (3.05,-1.35)
    {Положительная масса малого шара:\\[4pt]
     \(\nu(D_w)\ge m\,|B_{r/4}(0)|>0\).};
\end{tikzpicture}
\caption{Геометрия оценки роста строки. Цветом выделена часть
\(\overline B_r(y_c)\), лежащая по сторону \(\ip{e}{y-y_c}\ge r/4\)
от пунктирной границы. Весь малый шар лежит по эту сторону, а равенство
\(\norm{c_w-y_c}+r/4=3r/4<r\) обеспечивает \(D_w\subset B_r(y_c)\).
Схема двумерная; аргумент действует в \(\R^d\).
Исходный выбор \(\overline B_{2r}(y_c)\subset I\) на рисунке не показан.}
\label{ru:fig:bass-caps}
\end{figure}

Теперь нам нужно извлечь из роста строки положительное количество
кривизны. Нам не требуется, чтобы гессиан был положительно определён
в каждой отдельной точке: достаточно его массы на одном общем шаре.

\begin{lemmaru}\label{ru:lem:curv:hessian-mass}
Пусть выпуклая \(C^2\)-функция \(v\colon\R^d\to\R\)
\(M\)-липшицева, \(v(0)\le0\), и
\[
 v(w)\ge\ip{w}{y_c}+c\norm{w}-K,\qquad c>0.
\]
Тогда существуют \(T<\infty\), \(\eta>0\), зависящие только от
\(M,y_c,c,K,d\), для которых
\begin{equation}\label{ru:eq:curv:hessian-mass}
 \int_{B_T(0)}D^2v(w)[q,q]\,dw\ge\eta\norm{q}^2
 \qquad(q\in\R^d).
\end{equation}
\end{lemmaru}
\begin{proof}
Зафиксируем \(\rho>0\). При \(\norm{z}<\rho\), \(\norm{\xi}=1\)
из двух нижних оценок и \(v(z)\le M\norm{z}\) следует
\[
 v(z+L\xi)+v(z-L\xi)-2v(z)
 \ge2cL-2K-2(\norm{y_c}+M)\rho .
\]
Выберем \(L>0\) настолько большим, чтобы эта величина была не меньше
\(cL/2\), одновременно для всех таких \(z,\xi\).
Для \(j(t)=v(z+t\xi)\) выпуклость и основная теорема анализа дают
\[
 j(L)+j(-L)-2j(0)
 \le L[j'(L)-j'(-L)]
 =L\int_{-L}^{L}D^2v(z+t\xi)[\xi,\xi]\,dt.
\]
Поэтому интеграл справа без множителя \(L\) не меньше \(c/2\).
Проинтегрируем это неравенство по \(z\in B_\rho(0)\).
Подынтегральная функция неотрицательна. При каждом фиксированном \(t\)
её интеграл по переведённому шару не превосходит интеграла по
\(B_T(0)\), если \(T=L+\rho+1\). Теорема Фубини даёт
\[
 \frac c2|B_\rho|\le
 2L\int_{B_T(0)}D^2v(w)[\xi,\xi]\,dw.
\]
Можно взять \(\eta=c|B_\rho|/(4L)>0\).
Для произвольного ненулевого \(q\) применим оценку к
\(\xi=q/\norm{q}\) и воспользуемся билинейностью гессиана;
при \(q=0\) утверждение очевидно.
\end{proof}

По лемме~\ref{ru:lem:curv:cap-growth} одни и те же \(T,\eta\) подходят
всем интерполированным строкам при \(0<\varepsilon\le1\).
Именно независимость от \(\varepsilon\) позволит вернуться к исходной задаче.

\subsection{Гауссовское наблюдение и две копии}

Вернёмся к сравнивающей функции $J$ из раздела~\ref{ru:sec:ent:scalar}.
Здесь $Z_s=(1-s)Z_0+sZ_1$, $h=Z_1-Z_0$, а $b=a_1-a_0$.

Как отделить неизвестное приращение столбца \(b\) от приращения
исходного вектора \(h\)? Посмотрим на два независимых значения \(h\)
при одном и том же наблюдении. При вычитании выражений \(b\) исчезнет.

\begin{lemmaru}\label{ru:lem:curv:posterior}
Для \(T,\eta\) из леммы~\ref{ru:lem:curv:hessian-mass} положим
\[
 g_T=(2\pi)^{-d/2}e^{-(R+T)^2/2},
 \qquad b_\gamma=(2\pi)^{-d/2}.
\]
Тогда для всех \(0\le s\le1\)
\begin{equation}\label{ru:eq:curv:comparison-lower}
 J''(s)\ge\delta_R\,\E\norm{h-\E h}^2,
 \qquad \delta_R=\frac{g_T^2\eta}{2b_\gamma}>0 .
\end{equation}
\end{lemmaru}
\begin{proof}
Обозначим через \(\lambda_s=\Law(Z_s,h)\) закон пары на
\(\R^d\times\R^d\), а через \(\varphi\) --- стандартную гауссовскую
плотность. Поскольку оба \(Z_0,Z_1\) измеримы относительно \(\G\),
пара \((Z_s,h)\) независима от \(\Gamma\).
Закон наблюдения \(w=Z_s+\Gamma\) имеет плотность
\[
 \beta_s(w)=\int\varphi(w-z)\,\lambda_s(dz,dh)>0.
\]
Его условный закон скрытой пары задаётся явной вероятностью
\[
 \lambda_{s,w}(dz,dh)=
 \frac{\varphi(w-z)}{\beta_s(w)}\,\lambda_s(dz,dh).
\]
Формула проверяется интегрированием измеримых тестовых функций.
Мы строим её на конечномерном пространстве наблюдаемых величин;
регулярного условного распределения на самом исходном \(\Omega\)
здесь не требуется.

Зафиксируем \(w\) и положим \(p=p_{a_s,w}\).
В гильбертовом пространстве центрированных функций \(L^2(p)\)
обозначим через \(U_h\) центрирование функции
\(\ip{h}{y}-b(y)\). Для двух независимых условных копий \(h,h'\)
\[
 \norm{U_h-U_{h'}}_{L^2(p)}^2
 \le2\norm{U_h}_{L^2(p)}^2+2\norm{U_{h'}}_{L^2(p)}^2 .
\]
После усреднения получаем
\[
 \E_{\lambda_{s,w}}\operatorname{Var}_p(\ip{h}{y}-b(y))
 \ge\frac14\E_{\lambda_{s,w}^{\otimes2}}
       \operatorname{Var}_p\ip{h-h'}{y}.
\]
Столбец \(b\) сократился. По~\eqref{ru:eq:ent:hessian} последняя
дисперсия, делённая на \(\varepsilon\), равна
\(D^2v_{a_s}(w)[h-h',h-h']\).

На шаре \(B_T(0)\), поскольку \(\norm{z}\le R\), имеем
\[
 g_T\le\varphi(w-z)\le b_\gamma,\qquad
 g_T\le\beta_s(w)\le b_\gamma.
\]
После умножения двухкопийного условного интеграла на плотность
наблюдения его множитель становится
\[
 \frac{\varphi(w-z)\varphi(w-z')}{\beta_s(w)}
 \ge\frac{g_T^2}{b_\gamma}.
\]
Ограничим неотрицательный интеграл в
\eqref{ru:eq:curv:comparison-second} шаром \(B_T(0)\),
применим эту оценку и переставим интегралы:
\[
 J''(s)\ge\frac{g_T^2}{4b_\gamma}
 \int\!\!\int\!\!\int_{B_T(0)}
 D^2v_{a_s}(w)[h-h',h-h']\,dw\,
 \lambda_s(dz,dh)\lambda_s(dz',dh').
\]
Неотрицательность позволяет применять Тонелли; конечность нужных
интегралов следует также из непрерывности гессиана на замыкании шара
и ограниченности \(h\).
По~\eqref{ru:eq:curv:hessian-mass} правая часть не меньше
\[
 \frac{g_T^2\eta}{4b_\gamma}\E\norm{h-h'}^2
 =\frac{g_T^2\eta}{2b_\gamma}\E\norm{h-\E h}^2 .
\]
В последнем равенстве использовано обычное тождество дисперсии
для двух независимых копий вектора.
\end{proof}

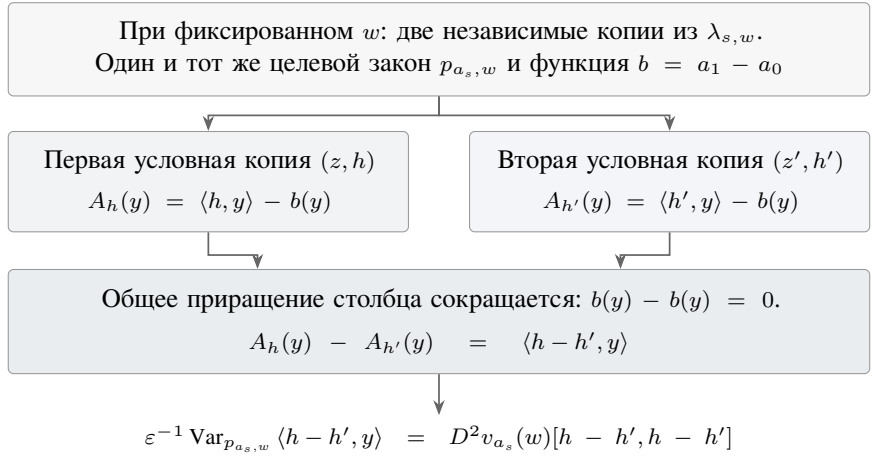
\begin{figure}[H]
\centering
\begin{tikzpicture}[
  font=\small,
  bass box/.style={draw=black!40,rounded corners=2pt,align=center,
    inner sep=7pt,minimum height=1cm},
  bass arrow/.style={-{Stealth[length=2mm]},draw=black!60,line width=.65pt}
]
\node[bass box,text width=10.9cm,fill=black!3] (obs) at (0,2.5)
  {При фиксированном \(w\): две независимые копии из \(\lambda_{s,w}\).\\
   Один и тот же целевой закон \(p_{a_s,w}\) и функция \(b=a_1-a_0\)};
\node[bass box,text width=4.8cm,fill=geoPtwo!6] (left) at (-3.05,.75)
  {Первая условная копия \((z,h)\)\\[3pt]
   \(A_h(y)=\ip{h}{y}-b(y)\)};
\node[bass box,text width=4.8cm,fill=geoPthreehalf!6] (right) at (3.05,.75)
  {Вторая условная копия \((z',h')\)\\[3pt]
   \(A_{h'}(y)=\ip{h'}{y}-b(y)\)};
\node[bass box,text width=10.9cm,fill=geoPtwo!9] (diff) at (0,-1.1)
  {Общее приращение столбца сокращается: \(b(y)-b(y)=0\).\\[4pt]
   \(A_h(y)-A_{h'}(y)=\ip{h-h'}{y}\)};
\node[align=center,text width=11.3cm] (hess) at (0,-2.65)
  {\(\varepsilon^{-1}\operatorname{Var}_{p_{a_s,w}}
       \ip{h-h'}{y}
     =D^2v_{a_s}(w)[h-h',h-h']\)};
\draw[bass arrow] (obs.south) -- ++(0,-.25) -| (left.north);
\draw[bass arrow] (obs.south) -- ++(0,-.25) -| (right.north);
\draw[bass arrow] (left.south) -- ++(0,-.3) -| ([xshift=-2.4cm]diff.north);
\draw[bass arrow] (right.south) -- ++(0,-.3) -| ([xshift=2.4cm]diff.north);
\draw[bass arrow] (diff.south) -- (hess.north);
\end{tikzpicture}
\caption{Что устраняет приём двух копий. При фиксированном наблюдении
\(w\) берутся две независимые копии евклидовой пары по одному
условному закону \(\lambda_{s,w}\).
Разность целевых потенциалов одинакова в двух выражениях и сокращается.
Дисперсия оставшегося линейного выражения связана с
гессианом логарифмического интеграла.}
\label{ru:fig:bass-two-copies}
\end{figure}

\subsection{Возвращение к функционалу Басса}

\begin{theoremru}\label{ru:thm:curv:bounded-slice}
Для каждого \(0\le R<\infty\) существует \(m_R>0\) такое, что при
\(Z,W\in\cH\), \(\norm{Z}_\infty,\norm{W}_\infty\le R\) и
\(\E Z=\E W\)
\begin{equation}\label{ru:eq:curv:strong-support}
 V(W)\ge V(Z)+\ip{g(Z)}{W-Z}_{\cH}
             +\frac{m_R}{2}\norm{W-Z}_2^2 .
\end{equation}
Константа выбирается до \(Z,W\) и их общего среднего.
\end{theoremru}
\begin{proof}
Для двух концов \(Z_0,Z_1\) лемма~\ref{ru:lem:curv:posterior} даёт
\(J''(s)\ge\delta_R\E\norm{h-\E h}^2\).
Следовательно,
\[
 J(s)\le(1-s)J(0)+sJ(1)
 -\frac{\delta_R}{2}s(1-s)\E\norm{h-\E h}^2 .
\]
Это обычная выпуклость функции
\(J(s)-(\delta_R/2)s^2\E\norm{h-\E h}^2\).
Заменим \(J(s)\) снизу регуляризованным значением, а концевые
\(J(0),J(1)\) --- равными им значениями из
\eqref{ru:eq:ent:primal-equality}. Вычитание
\(\E\ip{Z_s}{X}\) не меняет кривизну, поскольку это аффинная
функция \(s\).

Константа \(\delta_R\) не зависит от \(0<\varepsilon\le1\).
По утверждению~\ref{ru:prop:ent:removal} при
\(\varepsilon\downarrow0\) получаем
\[
 V(Z_s)\le(1-s)V(Z_0)+sV(Z_1)
 -\frac{\delta_R}{2}s(1-s)\E\norm{h-\E h}^2 .
\]
Если средние концов совпадают, \(\E h=0\), и последнее ожидание
равно \(\norm{Z_1-Z_0}_2^2\).
Наконец, положим \(Z_0=Z\), \(Z_1=W\), перенесём \(V(Z)\),
разделим на \(s>0\) и устремим \(s\downarrow0\).
Производная Фреше функционала даёт
\eqref{ru:eq:curv:strong-support} с \(m_R=\delta_R\).
\end{proof}

\noindent\textit{Комментарий и интуиция.}
Во всём пространстве функционал не может быть сильно выпуклым:
общий сдвиг исходного вектора сохраняет его значение, поскольку
барицентры \(\mu\) и \(\nu\) совпадают. Поэтому в оценке естественно
возникает центрированное приращение \(h-\E h\).
Фиксация среднего убирает именно это направление вырождения.
Гауссовский шум затем позволяет сравнить любые две ограниченные
конфигурации: на одном общем шаре наблюдений их плотности имеют
положительную общую нижнюю границу. Геометрия целевой меры обеспечивает
на этом шаре положительную массу кривизны. Вместе эти два факта и дают
положительное число \(m_R\).

\section{От кривизны к скорости потока}\label{ru:sec:rates}

Мы получили количественную оценку для двух ограниченных состояний с одинаковым
средним. Почему теперь её можно применять ко всем временам сразу?
Из качественной теоремы, указанной в разделе~\ref{ru:sec:setting}, существует один
$R<\infty$, для которого $\norm{Z_t}_\infty\le R$ при каждом $t\ge0$.
Сильный $L^2$-предел $Z^\star$ лежит в том же замкнутом множестве.
Средние всех этих векторов равны $\bar z=\E Z_0$.
Следовательно, всюду ниже действует одна и та же константа $m_R>0$
из оценки сильной выпуклости, а не новая константа для каждого времени.

Напомним её форму:
\begin{equation}\label{ru:eq:rate-support}
 V(W)\ge V(Z)+\ip{g(Z)}{W-Z}_{\cH}
       +\frac{m_R}{2}\norm{W-Z}_2^2\qquad(Z,W\in S_{R,\bar z}).
\end{equation}
В точке $Z^\star$ градиент равен нулю: $g(Z^\star)=0$.
Это следует из дифференцируемости $V$ и минимальности $Z^\star$ на всём $\cH$.

\subsection{Расстояние до предела}

Поменяем местами $Z,W$ в~\eqref{ru:eq:rate-support} и сложим обе оценки.
Значения функционала сократятся, и останется сильная монотонность:
\[
 \ip{g(Z)-g(W)}{Z-W}_{\cH}\ge m_R\norm{Z-W}_2^2.
\]
Подставляя $Z=Z_t$ и $W=Z^\star$, получаем
\begin{align*}
 \frac{d}{dt}\norm{Z_t-Z^\star}_2^2
 &=-2\ip{g(Z_t)-g(Z^\star)}{Z_t-Z^\star}_{\cH}\\
 &\le-2m_R\norm{Z_t-Z^\star}_2^2.
\end{align*}
Производная существует, поскольку поток является $C^1$-кривой.
Умножение на $e^{2m_Rt}$ показывает, что
$e^{2m_Rt}\norm{Z_t-Z^\star}_2^2$ не возрастает. При $t=0$ её значение
равно квадрату начального расстояния. После извлечения корня имеем
\begin{equation}\label{ru:eq:strong-norm-rate}
 \norm{Z_t-Z^\star}_2\le e^{-m_Rt}\norm{Z_0-Z^\star}_2.
\end{equation}
Этот аргумент также включает случай $Z_0=Z^\star$: обе стороны тогда равны нулю.

Оценка также контролирует законы в квадратичном расстоянии Вассерштейна.
Совместный закон $(Z_t,Z^\star)$ является допустимой связкой их распределений,
поэтому по определению этого расстояния
\[
 W_2(\Law(Z_t),\Law(Z^\star))\le\norm{Z_t-Z^\star}_2.
\]
Таким образом, контроль исходной гильбертовой траектории даёт и
количественное приближение искомого распределения.

\subsection{Зазор значений функционала}

Убывание расстояния само по себе ещё не даёт нужную оценку значений с точным
начальным множителем. Поэтому отдельно рассмотрим
$F(t)=V(Z_t)-V(Z^\star)\ge0$.
В~\eqref{ru:eq:rate-support} берём $Z=Z_t$, $W=Z^\star$. Если
$r_t=\norm{Z_t-Z^\star}_2$, то неравенство Коши--Буняковского даёт
\[
 F(t)\le\norm{g(Z_t)}_2r_t-\frac{m_R}{2}r_t^2
 \le\frac{\norm{g(Z_t)}_2^2}{2m_R}.
\]
Последний переход --- выделение полного квадрата:
$ar-m_Rr^2/2=a^2/(2m_R)-m_R(r-a/m_R)^2/2$.
Мы не минимизируем за пределами сечения: в опорную оценку был подставлен
конкретный допустимый вектор $Z^\star$.

По правилу дифференцирования сложной функции и уравнению потока
\[
 F'(t)=\ip{g(Z_t)}{\dot Z_t}_{\cH}
       =-\norm{g(Z_t)}_2^2\le-2m_RF(t).
\]
Тем же интегрирующим множителем получаем более сильную оценку
\begin{equation}\label{ru:eq:strong-energy-rate}
 V(Z_t)-V(Z^\star)
 \le e^{-2m_Rt}\bigl[V(Z_0)-V(Z^\star)\bigr].
\end{equation}
Поскольку $V(Z^\star)=\inf_{\cH}V$, выбор $\kappa=m_R$ одновременно
даёт обе оценки~\eqref{ru:eq:main-norm} и~\eqref{ru:eq:main-energy}.
Все остальные выводы основной теоремы уже обеспечены качественным результатом
раздела~\ref{ru:sec:setting}. Тем самым её доказательство завершено.

\paragraph{Комментарий и интуиция.}
Качественная теория удерживает траекторию в одной области. На этой области
равномерная кривизна не позволяет двум состояниям с одинаковым средним
оставаться далеко друг от друга при малом изменении градиента.
Поэтому расстояние до минимизатора убывает с единым экспоненциальным показателем.
Для энергии мы дополнительно использовали связь её зазора с квадратом нормы
градиента. Так получаются два разных дифференциальных неравенства и два
точных начальных множителя.

Мы не утверждаем оптимальности $m_R$ или показателя скорости. Константы
в предыдущем разделе выбраны для доказательства положительности и
равномерности; их численные значения могут быть далеки от лучших.

\section{Актуальность и связь с другими работами}
\label{ru:sec:related}

\subsection{Какой вопрос решает полученная оценка}

В задачах мартингального оптимального транспорта требуется не просто связать
два закона, но сделать это с сохранением условного среднего. Вариационный
подход к мартингалам Басса приводит к задаче минимизации поднятого функционала.
Для её изучения естественно следовать непрерывному градиентному спуску.
Однако существование предела и скорость приближения к нему отвечают на разные
вопросы: первое говорит, куда движется поток, второе количественно связывает
оставшуюся ошибку с прошедшим временем.

В работе Backhoff-Veraguas, Pammer и Schachermayer~\citep{ru:bps2025} качественный результат
получен в произвольной размерности, а две экспоненциальные оценки ---
в размерности один. Продолжение последнего утверждения на большие размерности
сформулировано там как гипотеза. Теорема~\ref{ru:thm:main} даёт это продолжение
при тех же условиях на данные. Известную качественную часть мы используем
явно; количественный вклад состоит в многомерной сильной выпуклости
на ограниченных сечениях фиксированного среднего.

\subsection{Непрерывный поток и дискретные методы}

Hasenbichler, Joseph, Loeper, Obłój и Pammer~\citep{ru:martingaleSinkhorn}
изучают мартингальный алгоритм Синкхорна. Его итерации строят потенциал
Басса в произвольной размерности. Это близкая задача калибровки,
но последовательность нелинейных обновлений не совпадает с потоком
\eqref{ru:eq:flow}. Теорема о сходимости таких итераций не определяет
автоматически скорость решения нашего дифференциального уравнения.

Marini~\citep{ru:marini2026} рассматривает более общий класс $q_x$-мартингалов
Басса и соответствующие вариационные и итерационные конструкции.
Эта работа расширяет семейство изучаемых мартингальных связок.
Здесь же фиксированы гауссовское сглаживание, функционал
\eqref{ru:eq:literal-functional} и его непрерывное время.
Мы не переносим наши оценки на весь класс $q_x$ без отдельного доказательства.

\subsection{Положительная энтропия и её исчезновение}

Связь с энтропийными задачами для нас особенно существенна, поскольку
энтропия входит в само доказательство. Hasenbichler, Pammer и
Thonhauser~\citep{ru:hasenbichler2026} изучают мост Шрёдингера--Басса,
его вариационное описание и предельные режимы. Этот круг результатов
объясняет, почему регуляризованные задачи естественно возникают рядом
с задачей Басса. Но сходимость их значений или оптимизаторов сама по себе
не даёт скорости конкретной нерегуляризованной траектории.

Количественный анализ энтропийных мартингальных алгоритмов представлен,
в частности, работами Kazeykina--Ren--Wang~\citep{ru:kazeykina2026}
и Bifronte--Ley--Sylvestre~\citep{ru:bifronte2026}.
Первая устанавливает экспоненциальную сходимость мартингального алгоритма
Синкхорна в относительной энтропии. Вторая даёт линейную сходимость SISTA
в произвольной размерности, сочетая контроль итераций, выпуклость
с учётом аффинных симметрий и оценку типа Поляка--Лоясиевича.
Эти работы важны для понимания общего механизма
«ограниченность плюс количественная кривизна».

Различие нашего аргумента находится в переходе $\varepsilon\downarrow0$.
Оценка второй производной сравнивающей функции имеет одну положительную
нижнюю константу при $0<\varepsilon\le1$. Благодаря этому сильная выпуклость
передаётся исходному функционалу. Одной скорости для задачи с фиксированным
положительным параметром было бы недостаточно: её константа могла бы исчезнуть
в пределе. Равномерность здесь доказана через гауссовское условное распределение
и интегральную массу гессиана на одном фиксированном шаре.

\subsection{Почему достаточно интегральной кривизны}

Для многомерного переноса непосредственный анализ изменения карты Бренье
может приводить к вопросам о регулярности нелинейной транспортной задачи.
Наше доказательство использует другой объект: логарифм интеграла,
гессиан которого выражается через ковариацию. Приём с двумя копиями
устраняет приращение столбца. После этого достаточно интегральной нижней
оценки гессиана, а не положительной поточечной границы во всём пространстве.
Поэтому доказательство не требует добавлять гладкость плотности или границы цели.
Это описание применённого метода; невозможности иных аналитических подходов
мы не утверждаем.

\paragraph{Ограничения.}
Основная теорема относится к ограниченным начальным данным и компактной
выпуклой цели с плотностью, ограниченной сверху и отделённой от нуля
почти всюду на внутренности поддержки, при указанной
неприводимости. Качественная теория допускает более общие старты, но
полученные здесь экспоненциальные оценки для всех таких стартов не доказаны.
Мы также не устанавливаем единую скорость для всех пар мер, оптимальный
показатель, скорость дискретизации или вычислительную сложность калибровки.
Обзор литературы разграничивает эти результаты, но не претендует
на исчерпывающий охват.

\section{Формальная верификация и воспроизводимость}\label{ru:sec:lean}

Итоговое многомерное утверждение формализовано в Lean~4 с использованием
mathlib. Проверка относится к точным исходным данным, функционалу и
дифференциальному уравнению, а не к абстрактному утверждению вида
«из сильной выпуклости следует сходимость».
В частности, формальная цепочка выводит количественную выпуклость,
коэрцитивность, радиальную оценку и ограниченность потока из исходных условий.
Готовый минимизатор или заранее ограниченная орбита в конечную теорему
дополнительными предпосылками не передаются.

Основной формальный результат носит имя
\begin{center}\small
\texttt{BassFlow.EntropicRoute.original\_Bass\_G}.
\end{center}
Вход \texttt{actual\_sourceBass\_G} задаёт тот же результат через данные мер.
В нём сохраняются исходное вероятностное пространство, закон $X$, независимый
гауссовский шум и произвольное ограниченное $\sigma(X)$-измеримое начальное
состояние. Вывод включает обе оценки с точными начальными множителями,
минимальность по всему исходному $L^2$ и единственность глобального минимизатора
при заданном среднем.

Полная проверяемая цепочка содержит 233 проектных модуля.
Исходники и инструкции находятся в самостоятельном репозитории
\url{https://github.com/Pure-Maths/Bass-Flow-Verification}.
С этой редакцией статьи сопоставлена версия доказательства, закреплённая
следующим коммитом:
\begin{center}\small
\href{https://github.com/Pure-Maths/Bass-Flow-Verification/tree/c8c9dec30e36585963df3f4a69508d30c26ff51f}{\texttt{c8c9dec30e36585963df3f4a69508d30c26ff51f}}.
\end{center}
Снимок этой закреплённой версии также приложен как дополнительный код.
Девять внешних модулей восстанавливаются описанным ниже шагом подготовки.

Файлы в каталоге \texttt{Numbered} раскрывают проверенные формальные
утверждения и их аксиоматические зависимости для сопоставления с нумерацией статьи.
Версии зависимостей зафиксированы в файлах \texttt{lean-toolchain}
и \texttt{lake-manifest.json}.

Для подготовки и проверки свежей копии используются команды
\begin{quote}\small\ttfamily
python3 scripts/bootstrap\_external.py\\
lake exe cache get\\
python3 scripts/check\_sources.py\\
python3 reproduce.py
\end{quote}
Первый шаг восстанавливает девять внешних модулей из зафиксированного
источника и проверяет контрольные суммы до и после совместимого переноса API.
В самом репозитории находятся остальные 224 проектных модуля.
Основные точки входа связываются с обычным доказательством следующим образом:

\begin{center}
\small
\renewcommand{\arraystretch}{1.2}
\begin{tabularx}{\linewidth}{@{}>{\raggedright\arraybackslash}p{.36\linewidth}>{\raggedright\arraybackslash}X@{}}
\toprule
Математический шаг & Основные модули формализации\\
\midrule
Градиент Фреше исходного $V$ &
\texttt{ActualSelectedContinuity}\\
Энтропийные оптимизаторы на концах и предельное неравенство &
\texttt{ActualEntropyChord} и его зависимости\\
Сильная выпуклость ограниченного сечения &
\texttt{ActualBassStrongConvexity}\\
Независимый вывод ограниченности и стационарного предела &
\texttt{SourceCenteredL1Coercivity},
\texttt{ActualBoundedSourceRadial},
\texttt{ActualBassStationaryFromRadial}\\
Минимальность, единственность и полный итог &
\texttt{ActualBassAllSourceUniqueness}, \texttt{ActualBassG}\\
\bottomrule
\end{tabularx}
\end{center}

В статье качественная часть явно заимствована из опубликованной
теоремы~1.8~\citep{ru:bps2025}. Формальный проект имеет более длинную
самостоятельную цепочку: он проверяет также построение потока, его удержание
в ограниченном множестве и минимальность предела на всём $\cH$.
Это различие порядка изложения не меняет итогового утверждения.
При сопоставлении с кодом необходимо смотреть формальные типы и зависимости;
таблица служит навигацией, а не заменой такого сопоставления.

Исходный пакет прошёл два раздельных аудита проверяемой цепочки.
В их отчётах для итоговых деклараций указаны только стандартные аксиомы
Lean: \texttt{propext}, \texttt{Classical.choice}, \texttt{Quot.sound}.
Проверка контрольных сумм, повторная компиляция проектных модулей
и доверие к зафиксированной внешней библиотеке являются разными уровнями
воспроизводимости. Команды и сведения о конкретных запусках приведены
в README самостоятельного проекта.

Формальная проверка не устанавливает новизну результата и не заменяет
научного рецензирования. Она также не делает автоматически формально
проверенной каждую поясняющую фразу текста или каждый рисунок.

\section{Использование ИИ при подготовке работы}\label{ru:sec:ai}

При работе над доказательством и его формальным представлением использовались
языковые модели. Они помогали прорабатывать технические переходы,
сопоставлять варианты аргумента и готовить Lean-код. При подготовке настоящей
русской редакции ИИ также использовался для построения структуры текста,
поиска и сопоставления первичных источников, подготовки LaTeX-исходников
и иллюстраций, а также независимых проверок отдельных частей доказательства.

При подготовке английской версии ИИ также использовался для перевода
и последующей сверки с русским текстом.

Для редакционной проверки привлекались несколько агентов с разными задачами:
математическая проверка, соответствие формальным утверждениям,
последовательность объяснений и техническая воспроизводимость.
Такие проверки помогают находить ошибки и неясности, но не заменяют
научного рецензирования. Проверка ядром Lean относится к точным формальным
формулировкам и их зависимостям; соответствие этих формулировок обычному
математическому тексту требует отдельного сопоставления.

ИИ не указывается в числе авторов. Выбор окончательных формулировок,
оценка математического содержания и ответственность за опубликованную
редакцию принадлежат автору.

\section{Заключение}\label{ru:sec:conclusion}

В этой работе мы рассматривали скорость сходимости непрерывного
градиентного потока поднятого функционала Басса.
При исходных условиях неприводимости, компактной выпуклой целевой
поддержке, плотности, почти всюду ограниченной сверху и отделённой от нуля
на внутренности поддержки, и ограниченном начальном состоянии получены
многомерные экспоненциальные оценки
расстояния в \(L^2\) и зазора функционала.
Одна положительная константа подходит для всех неотрицательных времён
данной траектории. Множители в оценках равны соответствующим
начальным ошибкам.

Основной количественный шаг состоит в сильной выпуклости на множествах
с единой верхней оценкой $L^\infty$-нормы и фиксированным средним.
Энтропийное сравнение позволяет записать кривизну через дисперсию.
Сравнение двух независимых копий устраняет изменение целевого
потенциала, а интегральная оценка гессиана даёт положительную нижнюю
границу, сохраняющуюся при удалении энтропии.
Известная качественная теорема помещает всю орбиту в одно ограниченное
множество. Это позволяет использовать одну константу кривизны на
всём промежутке времени и получить требуемую скорость.

Минимум рассматривается по всему исходному пространству случайных
векторов, а единственность относится к глобальным минимизаторам с данным
средним. Эти две оговорки имеют разные роли: первая задаёт область
сравнения, вторая учитывает симметрию функционала.

Также отметим границы результата. Мы не определяли оптимальный
показатель скорости и не доказывали существование одной скорости,
подходящей для всех допустимых мер и начальных данных.
Ограниченность начального данного и геометрические условия на
целевую меру остаются предпосылками теоремы.
При этом гладкость плотности и границы её поддержки не требуется.

\end{document}